\documentstyle[12pt]{article}
\begin{document}

\begin{center}

{\bf  BIDE - SIDE EXPONENTIAL AND }\\
\vspace{3mm}
{\bf MOMENT  INEQUALITIES FOR TAILS OF }\\
\vspace{3mm}

{\bf DISTRIBUTIONS OF POLYNOMIAL MARTINGALES.}\\

\vspace{3mm}

 {\sc BY OSTROVSKY E.I.} \\
\vspace{3mm}

  {\it  Ben - Gurion University, Beer - Sheva }\\

 \end{center}
\vspace{3mm}

  {\footnotesize In this paper non-asymptotic exponential estimates are derived for the 
  tail distribution of polynomial  martingale differences in  terms  unconditional 
  tails distributions of summands. Applications are considered in the theory of 
  polynomials on independent random variables, to
  the theory of $ U \ - $ statistics,
  multiply martingale series and in the theory of weak compactness measures on the 
  Banach spaces. \footnote{ Partially supported by the Israel Ministry of Absorbtion} 
\footnote{ Mathematics  Subject Classification (2000): 47A45, 47B10, 60F10, 60G42.} \\

 \normalsize

\vspace{2mm}

{\bf 1. Introduction. Notations. Statement of problem.} \ 
 Let $ (\Omega,F,{\bf P} ) $ be a probability space,
$ \xi(i,1), \xi(i,2), \ldots,\xi(i,d) $ be a family of centered 
$ ({\bf E} \xi(i,m) = 0) $ martingale - differences on the basis of the 
same flow  $ \sigma - $ fields (filtration) $ F(i): F(0) = \{\emptyset, \Omega \}, \ 
F(i) \subset F(i+1) \subset F, \ \forall m = 1,2,\ldots,d \Rightarrow
\xi(0,m) = 0; {\bf E} |\xi(i,m)| < \infty, $ and for every $ i \ge 0, 
m = 1,2,\ldots,d, \ \forall \ k = 0,1,\ldots,i-1 \ \Rightarrow  $

$$
{\bf E} \xi(i,m)/F(k) = 0; \ \ {\bf E} \xi(i,m)/F(i) = \xi(i,m) \ (mod \ \ {\bf P}),
$$

$ I = I(d) = \{(i_1,i_2,\ldots,i_d)\}, \ I(d,n) \ - $ {\it the set} of 
indexes $ I $ of the form $ I(d,n) = \{(i_1,i_2,\ldots,i_d)\}: $  
$ 1 \le i_1 < i_2 \ldots < i_{d-1} < i_d \le n, $ 
$ J(d) = J(d,n) $ - the subset of $ I(d,n) - $ the set of indexes of the form 
$ J(d,n) = J(d) = \{(i_1,i_2,\ldots,i_{d-1},n) \} $ \ such that $ 1 \le i_1 < i_2 
\ldots < i_{d-1} \le n-1, \ \ b(I) = b(i_1,i_2,\ldots,i_d) $ is a $ \ \  d \ - $ dimensional 
numerical non-random sequence,
$$
\xi(I)  = \prod_{m=1}^d \xi(i_m,m), \ \xi(J) = \prod_{m=1}^{d-1} \xi(i_m,m), \
\sigma^2(i,m) = {\bf D} \xi(i,m),
$$

$$
Q_d = Q(d,n, \{\xi(\cdot,\cdot)\}) = Q(d,n)  = \sum_{I \in I(d,n)} b(I) \ \xi(I) \ - \eqno(1.0)
$$
be a homogeneous polynomial (random polynomial) power $ d $ on the variables 
$ \{\xi(i,m) \} $ "without diagonal members",  (on the other hand, multiply 
stochastic integral on the discrete martingale measure, martingale transform), 
$ n $ is an integer number: $ n = 1,2,\ldots, \infty; $ in the case $ n = \infty $  
$ Q(d,\infty) $ should be understood 
 as a limit \ \ $ Q(d,\infty) = \lim_{n \to \infty} Q(d,n) $ (with probability 1). \par
 Note than 

$$
{\bf D} Q(d,n) = \sum_{I \subset I(d,n)} b^2(I) \prod_{m=1}^d \sigma^2(i_m,m);
$$
hence, if
$$
\sum_{I \subset I(d,\infty)} b^2(I) \prod_{m=1}^d \sigma^2(i_m,m) < \infty, \eqno(1.1)
$$
then by virtue of theorem of D.Doob $ Q(d,\infty) $  exists. In particular, if 
$$
\sup_{i,m} {\bf D} \xi(i,m) = \sup_{i,m} \sigma^2(i,m) < \infty, \
\sum_{I \subset I(d,\infty)} b^2(I) < \infty, \eqno(1.2)
$$
then condition (1.1) is satisfyed.\par
  We shall denote 
$$
B(d,n)= B = \{b(I): \sum_{I \in I(d,n)} b^2(I) = 1 \}, \ \ n \le \infty \eqno(1.3) 
$$
and assume that $ b \in B = B(d,n). $ \par
 If $ b(\cdot) \in B, \ \sigma^2(i,m) = 1, $ then $ {\it D} Q_d = 1. $ \par 
Another notations.  For any random variable $ \tau $ we define
$$
 T(\tau,x) = \max \left({\bf P}(\tau > x), {\bf P}(\tau < -x) \right), \ x > 0, \ - 
$$
 the tail of distribution $ \tau, $

$$
T_{Q(d,n)}(x) = T_Q(x) = T(Q_d,x).
$$
 {\bf Our goal is to make  non - asymptotic uniform over $ b \in B $ 
estimate of the tails of random variables
 $ T(Q_d,x) $ and of the moments $ {\bf E}|Q_d|^p $  in the terms of unconditional tails  
$ T(\xi(i,m),x) $ and moments  $ {\bf E}|\xi(i,m)|^p $ of  summands } $ \{\xi(i,m) \}.$ \par
 To interest our readers, we shall  formulate now two simple results. 
 Denote as $ G(q), \ q > 0 $ the set of all random variables $ \{\eta\}, $ 
which are defined on the our probability space $ \{\Omega,F,{\bf P} \}, $ such that 
$$
\exists K = const \in (0,\infty], \ \forall x > 0 \ \Rightarrow  \
T(\eta,x) \le \exp \left( -(x/K)^q \right).
$$
 It follows  from  theory of $ G \ - $ spaces ([11], p. 31 - 37) that in according to the norm 
$$
||\eta||_q := \sup_{m \ge 1} \left[ \ |\eta|_m \ \right] \ m^{-1/q}, \ \ |\eta|_m 
\stackrel{def}{=} \left[{\bf E}|\eta|^m \right]^{1/m}
$$
the set $ G(q) $ is the (full) Banach space which is isomorphic to the Orlicz 
space on the probability space $ (\Omega,F, {\bf P}) $ with $ N \ - $ function $ N(u) = 
\exp(|u|^q) - 1 $ [11, p. 35 - 37]. \par
 In the case $ q = \infty $ the space $ G (\infty) $ consists of all bounded 
$ (mod \ {\bf P} ) $ variables, and the norm $ G(\infty) $ is equivalent to the 
classical $ L_{\infty} $ norm
$$
|\eta|_{\infty} = vraisup_{\omega \in \Omega}|\eta(\omega)|.
$$

 Generalization: denote $ G(q,r), \ q > 0, r \in (-\infty,+\infty) $ the set of all random 
variables $ \{\eta\} $ with finite norm:
$$
||\eta||_{q,r} = \sup_{p \ge 2} |\eta|_p  \ p^{-1/q} \ \log^{-r}p; 
$$
it is known that $ ||\eta||_{q,r} < \infty $ if and only if 
$$
   \exists K = const >0, \  T(\eta,x) \le \exp  \left(-(x/K)^q(\log(F + x/K))^{-qr}\right),
$$
and for some $ C_1,C_2 = C_1(q,r), \ C_2(q,r) \in (0,\infty), \ C_1 \le C_2 $
$$
  C_1(q,r) \ K \le ||\eta||_{q,r} \le C_2(q,r) \ K, 
$$
where   
$$
F = F(q,r)=1, \ r \le 0;  \ \ F(q,r)= \exp(q), \ r > 0.
$$
 {\it Throughout this  paper the letter $ C_j(\cdot) $ will denote various constants which 
may differ  from one formula to the next even within a single string of estimates and which 
do not depend upon $ n.$ We make no attempt to obtain 
the best values for these constants.}\par 
 Another example (concise). Let us introduce the other space of random variables 
$ \Psi(C,\beta), \ \beta = const > 0. $  By definition, this space
 consist of all random variables $ \{\eta\} $ with finite norm
$$
|||\eta|||_{C,\beta} \stackrel{def}{=} \sup_{p \ge 1} |\eta|_p \exp \left(-Cp^{\beta} \right); 
\ \eta  \in \cup_{C>0} \Psi(C,\beta) \ \Leftrightarrow 
$$
$$
 \exists C_1 > 0, \  T(\eta,x) \le \exp \left(-C_1 (\log (1+x))^{1+1/\beta} \right).
$$

 For instance assume that there exist constants $ K(m) \in (0,\infty), \
q(m) > 0 $ so that $ \forall x \ge 0  $

$$
\sup_i T(\xi(i,m),x) \le \exp \left( - (x/K(m))^{q(m)} \right), \eqno(1.4)
$$
or  briefly $ \max_{m=1,2,...,d} \sup_{i \le n} ||\xi(i,m)||_{q(m)} /K(m) < \infty.$
 We define $ K = \prod_{m=1}^d K(m), \ 1/\infty = 0, \ M = M(d,q)= $
 
$$
 M(d,\vec{q}) = M(d;q(1),q(2),\ldots,q(d)) = 
\left( d/2 + \sum_{m=1}^d (1/q(m)) \right)^{-1},
$$
$$
S(\vec{q}, \ x) = S(q(1),q(2),\ldots,q(d);x)= \sup_{b \in B} \sup_{\{\xi(i,m)\}} T_Q(x),
$$
where interior  $ {\bf \sup} $ is calculated over all families of centered  martingale 
- differences $\{\xi(i,m)\} $ satisfying  condition (1.4).\par

{\bf Theorem 1.} {\it  There exist two constants} $ C_1 = C_1(\vec{q}), \ C_2 =
C_2(\vec{q}), \ 0 < C_2 < C_1 <\infty $ {\it so that for  all } $ x > x_0 = const > 0 $

$$
\exp \left(- \left[x/(C_1 K) \right]^M  \right) \le S(\vec{q},x) \le \exp 
\left(- \left[x/(C_2 K)^M \right] \right). \eqno(1.5).
$$
 On the other hand, the right inequality in our result (1.5) may be rewritten in 
the terms of $ G(q) \ - $ spaces as
$$
\sup_{ \{\xi(i,m): ||\xi(i,m)||_{q(m)}< \infty \} }  ||Q_d||_{M(d,\vec{q})} /
\prod_{m=1}^d \sup_{i \le n} ||\xi(i,m)||_{q(m)}   \le C_3(\vec{q}), \eqno(1.6)
$$
and the left part of Theorem 1 denotes that inequality (1.6) can not be improved.\par
 We can improve the result (1.5) in so - called " independent case", i.e. if all the
 variables $ \{\xi(i,m) \} $ are independent.  Let us  consider a more 
generally centered polynomial (with "diagonal members") degree $ d $ of a view 
$ R_d=R_d(n) = R_d(n, \{\xi(i,s) \}) = $ 
$$
 \sum_{1 \le i_1 < i_2 \ldots < i_d \le n}
b(i_1,i_1,\ldots,i_1;i_2,i_2,\ldots,i_2; \ldots,i_d,i_d,\ldots,i_d) \times
$$
$$
 \prod_{l=1}^d \left(\xi^{k(l)}(i_l,l) - m(k(l),i_l,l) \right),
$$
where  $ m(k,i,l) = {\bf E} \xi^k(i,l), \ k(l) \in \{0,1,\ldots,d\},$
so that $ \sum_l k(l)  = d, $
$$
k(l) = card \{i_l \ {\bf in} \ \ b(i_1,i_1,...,i_1;..., i_l,i_l,...,i_l; ...,
i_d,i_d,...,i_d) \},
$$
and if some $ k(l) = 0, $ then by definition,
$$
\xi^{k(l)}(i_l,l) - m(k(l),i_l,l) = 1.
$$
 We again suppose that $ b \in B;  $  (the sequence $ b(I) $ may be not  
symmetric.) The multiply series for $ R_d $ 
in  the case $ n = \infty $  converge with  probability 1, for instance,
if $ b \in B $ and the correspondence even moments are bounded:
$$
  \sup_{\{i,l\}} m(2d,i,l) < \infty. \eqno(1.7)
$$
 For example if $ d = 2 $ then $ R_d $ has a view 
$$
R_2 = \sum \sum_{1\le i < j \le n} b(i,j)\xi(i)\xi(j) + \sum_{i \le n} b(i,i)
(\xi^2(i) - {\bf E} \xi^2(i)) + 
$$
$$
\sum \sum_{ \le i < j \le n} c(i,j) \eta(i)\eta(j) + \sum_{i \le n} c(i,i)
(\eta^2(i) - {\bf E}\eta^2(i)) +
$$
$$
 \sum \sum_{1 \le i,j \le n} a(i,j) \xi(i)\eta(j).
$$
 If $ d = 3, $ then

$$
R_3 = \sum \sum \sum_{1 \le i < j < k \le n} b(i,j,k) \xi(i) \eta(j) \tau(k) +
$$
$$
 \sum \sum_{1 \le i < k \le n} b(i,i,k) (\xi^2(i) - {\bf E} \xi^2(i)) \tau(k) +
\sum_i b(i,i,i) (\xi^3(i) - {\bf E} \xi^3(i) ) + \ldots.
$$
 Here $ \{\xi(i), \eta(j), \tau(k) \} $ are independent sequences of 
independent random variables.\par
 Denote 
$$ 
U(\vec{q},x) = U(q(1),q(2),\ldots,q(d);x) = 
 \sup_{b \in B} \sup_{ \{\xi(i,m) \in G(q(m)) \} } T(R_d,x),
$$
where  $ \sup_{ \{\xi(i,m) \in G(q(m)) \} } $ is calculated over all totally 
independent random variables $ \{\xi(i,m) \} $ such that  $ \exists q(m) > 0, \ 
K(m) > 0 \ \Rightarrow \forall m = 1,2,\ldots,d $
$$
\sup_i T(\xi(i,m),x) \le \exp \left(-(x/K(m))^{q(m)} \right), \ x > 0,
$$
or, equally  

$$
\sup_i ||\xi(i,m)||_{q(m)} \le C \ K(m), \ C = const < \infty.
$$

 Put $ K = \prod_{m=1}^d K(m), \ N(q) = 2q/(q+2) $ by $ q \in (0,1], $
$ N(q) = \min(q,2) $ if $ q > 1,$  and define a family of function 
$ N(\vec{q}) = N_d(\vec{q}) = N_d(q(1),q(2),\ldots,q(d)) $ by the following 
recursion: $ N_1(q) = N(q), $ (initial condition),

$$
N_{d+1}^{(k)}(q(1),\ldots,q(d),q(d+1)) = \left[\frac{d-1}{2} + 
\sum_{m=1,2,...,d; m \ne k} \frac{1}{q(m)} + \frac{1}{N(q(k))} \right]^{-1},
$$
$$
N_{d+1}(q) =  N_{d+1}(q(1),q(2),\ldots,q(d), q(d+1)) =
$$
$$
 \max_{k=1,2,\ldots,d+1} N^{(k)}_{d+1}(q(1),q(2),\ldots,q(d+1)).
$$
{\bf Theorem 2.} {\it There exists a function} $ C_3 = C_3(\vec{q}) \in (0,\infty) $ {\it 
such that } $ \forall x \ge 2 $
$$
U(q(1),q(2),\ldots,q(d),x) \le \exp \left(- C_3(\vec{q}) \ (x/K)^{N_{d}(q)} 
\right).\eqno(1.8)
$$
 In the terms of $ G(q) \ - $ spaces the proposition (1.8) may be rewritten as 
$$
||R_q||_{N_d(q)} \le C_4(d,\vec{q}) \prod_{m=1}^d \sup_{i \le n}
||\xi(i,m)||_{q(m)}.
$$

 For example  assume that $ q(m) = q = const > 0 $ and denote 
$$
\gamma(d,q) = 2q/[d(q+2)], \ \ q \in (0,1],
$$
$$
\gamma(d,q) = 2q/[2d+q(d-1)],  \ \ q \in (1,2],
$$
$$
\gamma(d,q) = 2q /[d q + 2(d-1)], \ \ q \in (2,\infty].
$$
 In this case $ N_d(\vec{q}) = \gamma(d,q), $ and we receive  the following 
result: if $ \xi(i,m) $ are totally independent, centered  and 
$$
T(\xi(i,m),x) \le \exp \left(-x^q \right), \ x \ge 0,
$$
then by all $ x \ge x_0 $  and $ \ dim(q,q,\ldots,q) = d \ \Rightarrow  $ 
$$
U(q,q,\ldots,q,x) \le \exp \left(-C_3 \ x^{\gamma(d,q)} \right). \eqno(1.9)
$$

 {\bf Theorem 3.} (Low bounds for $ U(q,q,...,q,x) $ ).  {\it There exists a function} 
$ C_4 = C_4(d,q) \in (0,\infty) $ {\it such that for all} $ x \ge 1 $
$$
U(q,q,\ldots,q, x) \ge \exp \left(-C_4 \cdot x^{\min(q,2)/d} \right).\eqno(1.10)
$$
 Obviously upper estimations (1.9) and low (1.10) "almost" coincides: at $ q \to 
0+ $ or by $ q = \infty $ for all $ d, $  and by $ d = 1, \ q \in (1,\infty].$
If for example $ \{\xi(i,s) \} $ are independent Rademacher series: 
$ {\bf P}(\xi(i,s) = 1) = {\bf P}(\xi(i,s) = - 1) = 1/2 $ and $ \sum_I b^2(I)
< \infty, $ we deduce the well - known result  ([28, p. 78], [31]) as a particular case:
$$
\exists \varepsilon > 0; \ {\bf E} \exp \left(\varepsilon |Q_d|^{2/d} \right) 
< \infty.
$$
 Furher we shall  formulate and prove more generally results.\par
 There are many publications about the limit theorem, moment and exponential 
inequalities for tail distributions of martingales and polynomials from 
independent variables. Semiinvariant inequalities for random polynomials are received 
in the book [15, p. 100 -103] in the case Gaussian 
limit distribution.  The case $ d = 2 $ is considered 
in paper [30] and it is proved  that in all rearrangement 
invariant space $ X, $ for example  $ X = G(q), $ 
$$
||\sum \sum_{1 \le i \ne j \le n} b(i,j) \xi(i,1)\xi(j,1) ||_X \asymp
\sqrt{ \sum \sum b^2(i,j)}.
$$
 Consequently, the norm $ \sqrt{\sum \sum b^2(i,j)} $ is optimal. We can 
explain this from identity $ | \sum_i b(i) \xi(i)|_2^2  =
 {\bf D}\sum b(i)\xi(i) = \sum b^2(i), $ where 
$ \{\xi(i)\} $ are  sequence of centered martingale - differences 
 with condition $ {\bf D}\xi(i) = 1. $\par
 Non - uniform estimations $ T_Q $ are obtained 
in the works [1], [2], [3], [6], [7] etc.  They are received in the terms 
of {\it conditional} expectations $ {\bf E} g(\xi(i))/F(i-1) $ of a view, for 
example:
$$
{\bf P}(\exists n, \ \sum_{i=1}^n \xi(i) > x, \ \sum_{i=1}^n {\bf E} \xi^2(i)
/F(i-1) \le y ) \le 
$$
$$
\le \exp (- x^2/(2(y+Cx))), \ x,y > 0,
$$
i.e. without sequences of coefficients $ b(I) $ as in classical Bernstein - 
Bennet estimations. But these estimations are not convenient in the practice.\par
 So - called "decoupling method" for calculation of order of magnitude expectation
$$
{\bf E} f(\sum \sum_{i,j=1}^n g_{i,j}(\xi(i), \eta(j)) 
$$
is described in the  articles [2], [4] and in other publications,
but only for a function $ f, $ belonging  $ \Delta_2 \ - $
class. See another publications in references. \par
 The limit theorems for  martingales are well - known [20, p. 58]. 
 The limit theorems are received in [24] for symmetric polynomials, 
 for example of a view: 
$$
Q_d = \sum \sum \ldots \sum_{1 \le i_1 < i_2 \ldots < i_d \le n} \prod_{i=1}^d 
\xi(i_s),
$$
where $ \{\xi(i) \} $ are i.,i.d. random variables and it is proved that under some
conditions at $ n \to \infty $ in the sense of distribution convergence
$$
Q_d(n) /\sqrt{ {\bf D} Q_d(n)} \stackrel{d}{\to} \ I(h), \eqno(1.11)
$$
where $ I(h) $ is a multiply stochastic integral 
$$
I(h) = \int \int \ldots \int_{R^d} h(\lambda_1,\lambda_2,\ldots,\lambda_d)
\prod_{s=1}^d Z(d\lambda_s),
$$
$ Z(\cdot) \ - $ white Gaussian measure: $ {\bf E} Z(A)=0, \ {\bf E} Z(A) 
Z(B) =  mes (A \cap B), \ h \in L_2(R^d), h \ne 0. $ \par 
 {\it Our estimations are formulated in the very simple terms only {\it unconditional}
individual (marginal) tails of summand distributions, are very convenient for 
using, generalized in the multidimensional case $ d > 1, $ and are non - 
improved essentially (see, for example, theorem 1 (1.5) and inequalities (1.9), (1.10)).} \par

\vspace{3mm}
{\bf 2. Main results: exponential estimations.} \
  We shall  introduce some notation used in following sections.  Let  $ T(x) $ and $ G(x), 
x > 0 $ be two tail - functions, i.e. $ T(0) = G(0)=1, $  monotonically 
decreasing  right continuous and so that $ T(\infty) = G(\infty) = 0. $ We denote
$$
T \vee G(x) =   \min(4\inf_{y > 0} (T(y)  + G(x/y) ),1).
$$
 The function $ T \vee G(x) $ has the following sense: if $ T(\xi,x) \le T(x), \ 
T(\eta,x) \le G(x) $ then 
$$
 T(\xi \cdot \eta,x) \le T \vee G(x). \eqno(2.1)
$$
For example if $ \ \forall \ x \ge 0 \ T(\xi,x) \le Y_1 \exp \left( -(x/A)^{q(1)} \right), \
T(\eta,x) \le Y_2 \exp \left( -(x/B)^{q(2)} \right) $ for some $ q(1), q(2), A,B = const > 0, 
Y_1, Y_2 = const \ge 1, $ then $ \forall x \ge 0 \ \Rightarrow $ 
$$
T(\xi \cdot \eta, \ x) \le 8 \max(Y_1,Y_2) \ \exp \left( - (x/(AB))^{q(1) q(2)/(q(1)+q(2))}
\right).
$$

 More generally, if
$$
T(\xi,x) \le \exp \left(-x^{q(1)} (\log(F(q(1),r(1)) + x))^{r(1)} \right), 
$$
$$
T(\eta,x) \le \exp \left(-  x^{q(2)}(\log(F(q(2),r(2))+x))^{r(2)} \right),
$$
then  $ T(\xi \cdot \eta, x) \le $
$$
  \min(1,4 \exp \left(- C x^{q(3)} \ (\log(F(q(3),r(3)) + x))^{r(3)} \right),\eqno(2.2)
$$
 where
$$
q(3) = q(1)q(2)/[q(1)+q(2)], \ \ r(3) = [q(1)r(2) + q(2)r(1)]/[q(1)+q(2)].
$$
 On the other hand,  according to the
 language of $ G(q,r) \ - $ spaces: for some $ C = C(q(1),q(2),r(1),r(2)) \in 
(0,\infty) \ \Rightarrow $
$$
||\xi \cdot \eta||_{q(3),-r(3)/q(3)} \le C ||\xi||_{q(1),-r(1)/q(1)} \cdot ||\eta||_{q(2),-r(2)/q(2)}. 
$$
 If $ \xi,\eta $ are independent then the  estimation (2.2) is exact.\par 
 Further, let us denote for the tail - function $ T(\cdot) $  the following operator 
(non - linear)
$$
W[T](x) = \min \left(1, \inf_{v > 0} \left[ \exp(-x^2/(8v^2)) - \int_v^{\infty}
x^2 \ dT(x) \right] \right),
$$
 if there exists the second moment $ |\int_0^{\infty} x^2 \ dT(x)| < \infty. $ \par
{\bf Lemma 1. } {\it Let} $ d=1, \ \xi(i) = \xi(i,1) $  {\it be a sequence of martingale - 
differences with filtration} $ \{F(i)\} $ {\it and let } $ T(\xi(i),x) \le T(x),$ 
where $ \ T(x) $ {\it is some tail - function. Then for all} $ x \ge 2 $ 

$$ 
\sup_{b: \in B} T \left(\sum_i b(i) \xi(i), \ x \right)  \le W[T](x). \eqno(2.3)
$$
{\bf Proof}. Without loss of generality we can assume that $ i = 1,2,\ldots,n; $
where $ n < \infty. $
We shall use the so - called "truncation" method [10]. 
 The following inequality  for finite martingale differences is known, i.e. 
when $ vraimax_{\omega} |\xi(i)| = c(i) < \infty $  (see [8], [25]) 
$$
T(\sum \xi(i),x) \le \exp \left(-x^2/\left(2\sum c^2(i) \right) \right). \eqno(2.4)
$$
Note that in [25] the factor 2 is omitted in the denominator of the exponent 
index in the corresponding formulation. \par
 Put  $ X = \sum_i b(i) \xi(i), \ Y(1,i) = b(i) \xi(i) \chi(|\xi(i)| \le v), $
$$
Y(i) = Y(1,i) - {\bf E}Y(i,i)/F(i-1), 
$$
$$
Z(1,i)  = b(i) \xi(i) \chi(|\xi(i)| > v), \ v = v(x) > 0,
$$
$$
Z(i) = Z(1,i) = {\bf E}Z(1,i)/F(i-1), 
$$

$$
\chi(A) = 1, \ \omega \in A, \ \chi(A) = 0, \ \omega \notin A, \ A \in F.
$$
 We can write $ X = X_1 + X_2, $ where $ X_1 = \sum Y(i), X_2 = \sum Z(i), $ and we obtain 
according to the inequality (2.4) and since $ \ |Z(i)| \le |b(i)| v: $   
$$
T(X_1,x) \le \exp \left( - x^2/(2\sum b^2(i) v^2) \right) = \exp \left( - x^2/2v^2 \right).
$$
 Further, if $ b(i) \ne 0 $ then $ \ {\bf D}Z^2(i) b^{-2}(i) =  {\bf E} Z^2(i) 
b^{-2}(i)  \le $
$$
 - \int_v^{\infty} y^2 \ dT(\xi(i),y) = v^2 T(\xi(i),v)  + 2 \int_v^{\infty}
y T(\xi(i),y) \ dy \le 
$$
$$
\le v^2 T(v) + 2 \int_v^{\infty} yT(y) \ dy  = - \int_v^{\infty} y^2 \ dT(y).
$$
 Since $ \{Z(i)\} $ are not correlated, we see that
$$ 
{\bf E} X_2^2 \le - \sum_i b^2(i)  \int_v^{\infty} y^2 \ dT(y)  =
 - \int_v^{\infty} y^2 \ dT(y).
$$
By virtue of  Chebyshev inequality 

$$
T(X_2,x) \le \left| \int_v^{\infty} y^2 dT(y) \right| /x^2.
$$
 Now our statement it follows  from the simple inequality that 
$$
{\bf P}(X > x) \le {\bf P} (X_1 > x/2) + {\bf P}(X_2 > x/2), \ x \ge 2,
$$
after the minimization on  $ v.$  This completes the proof of lemma 1.\par

 For example, assume   $ T(\xi(i),x) \le Y \ \exp(-(x/K)^q), \ K,x,q \in (0,\infty),
 \ Y  = const \ge 1. $ Denote (in time, in this section) 
$$
\delta = \delta(q) = (\min(q/2, \ 1)^{-1/q}; \ \ q \in (0,2] \ \Rightarrow 
$$

$$
\beta(q) = \max \left( (1/q) \Gamma (2/q), \ (e/q) [2/(eq)]^{2/q} \right),
$$

$$
\beta(q) := \sup_{v \ge 0} \ \exp \left(v^q \right) \ \int_v^{\infty} x 
\exp \left(-x^q \right) \ dx \le \Gamma (2/q)/(qe)
$$
in the case $ q > 2; $ here $ \Gamma(\cdot) $ is the Gamma - function. After some calculation 
we find for all values $ x \ge 0:$
$$
\sup_{b \in B} T(Q_1,x) \le W[T](x) \le (1+2 \beta(q) \ Y) \ \exp \left( - (x/(K \delta))^{2q/(2+q)} 
\right).
$$
 Analogously assume that $ \sup_i ||\xi(i)||_{q,r} =K < \infty, q > 0, 
r \in (-\infty,\infty) $ or equally $ x \ge 2 \ \Rightarrow $ 
$$
T(\xi(i),x)  \le \exp \left(- (x/K)^q \ (\log(F(q,r) + x/K))^r \right).
$$
 Let us introduce the following vector - function $ L(q,r) = \{L(1;q,r), \ L(2;q,r)\}  $ 
of two variables $ \ (q,r):$
$$
L(1;q,r)= 2q/(q+2), \ \ L(2;q,r) = 2r/(q+2).
$$
 We see that for $ x > 0:\ \sup_{b \in B}  T(\sum_i b(i)\xi(i),x) \le $
$$
 \exp \left(- (x/(C_1K))^{L(1)} \ (\log(F(L(1),L(2))+x/(C_1 K))^{L(2)} \right), \eqno(2.5) 
$$
where $ L(i) = L(i;q,r), \ i = 1,2; $ or 
$$
\sup_{b \in B} ||\sum_i b(i)\xi(i)||_{L(1;q,r),L(2;q,r)} \le C_2(q,r) \sup_i ||\xi(i)||_{q,r}.
$$
 The proposition (2.5) is new even for independent variables 
$\{\xi(i)\} $ in the case $ q \in (0,1]. $ \par
{\bf Theorem 4.} (Martingale case). {\it Let us  denote } $ T_m(x) = \sup_i T(\xi(i,m),x) $ 
{\it and assume that} $ \lim_{x \to \infty} T_m(x) = 0, $ {\it and define the sequence} 
 $ T^{(s)} $ {\it of tail functions in the following way:} 
$$
T^{(1)} (x) = W[T_d](x), 
$$
{\it (initial condition), and for } $ s = 2,3,\ldots,d-1 $ {\it by recurrent 
equation }
$$
T^{(s+1)}(x) = W \left[T_{s+1} \vee T^{(s)}  \right](x).
$$
 {\it Statement:} 
$$
\sup_{b \in B} T(Q_d,x) \le T^{(d)}(x). \eqno(2.6)
$$
 {\bf Proof}. We shall prove the statement (2.6)  by means of induction over $ d. $  
 The basis of induction  $ (d = 1) $ in (2.6) is proved in  lemma 1. Further, 
the sequence $ \{ Q_d(n), F(n) \} $ is again martingale with correspondence 
martingale - differences 
$$ 
\zeta(n) = Q_d(n) - Q_d(n-1) = \xi(n,d) \sum_{I \in J(d,n)} b(I) \xi(I) =
$$

$$
= \xi(n,d) \sum_{1 \le i_1 < i_2 \ldots < i_{d-1} \le n-1 } b(i_1,i_2,\ldots, 
i_{d-1},n)\prod_{m=1}^{d-1} \xi(i_m,m) = 
$$

$$
= \left[\xi(n,d)\times \sqrt{\sum_{I \in J(d,n)} b^2(I) } \right] \times 
\left[\sum_{I \in J(d,n)} \frac{\sum_{I \in J(d,n)} \xi(J)b(I)}
{ \sqrt{\sum_{I \in J(d,n)}b^2(I) }   } \right] \stackrel{def}{=}
$$
$$
\stackrel{def}{=} [\eta(n,d)] \times [\tau(n,d)].
$$
 We get  in  according of induction statement:
$$
T(\tau(n,d),x) \le T^{(d-1)}(x).
$$
 We deduce by virtue of the formula (2.1) for the product $ \eta(d,n) \times \tau(d,n) $ 
$$
T(\eta(d,n)  \cdot  \tau(d,n),x) \le \left( T_d \vee T^{(d-1)} \right)(x).
$$
 Again using the statement of lemma 1 for a martingale differences 
$ \eta(d,n) \ \tau(d,n) $ we obtain the statement of theorem 4. \par
 If for example
$$ 
T(\xi(i,m),x) \le \exp \left( - (x/K(m))^{q(m)} \right), 
$$
 After some calculations we shall receive the statement of theorem 1.\par
{\it Let us consider now "independent" case.} Define
$$
\varphi_m(\lambda) = \sup_i \max_{\pm} \log {\bf E} \exp(\pm \lambda \xi(i,m)).
$$
 This definition is non - trivial:  
$$
 \exists \lambda_0  \in (0,\infty], \ \forall \lambda \in (-\lambda_0, \lambda_0) \ 
\Rightarrow \varphi_m(\lambda) < \infty \eqno(2.7) 
$$
only if the sequence $ \{\xi(i,m) \} $ satisfies the uniform Kramer condition. \par
 If the condition (2.7) is valid, we introduce the functions  $ \chi_m(\lambda) $ by formula
$$
\chi_m(\lambda) = \sup \left\{ \sum_{j=1}^{\infty} \varphi_m(\lambda y(j)), \
\{y(j)\}: \sum_{j=1}^{\infty} y^2(j) \le 1 \right\}.
$$
 Note that if the function $ z \to \varphi_m^/(z)/z $ is finite and monotonic 
on the right - hand half - line then 
$$
\chi_m(\lambda) = \sup_n n \ \varphi_m(\lambda/\sqrt{n}). \eqno(2.8)
$$
 Namely:
$$
\chi_m(\lambda) = \sup_n \max_{\{y(j)\}} \sum_{j=1}^n \varphi_m(\lambda y(j)),
$$
where the internal maximum is taken over all finite sequences $ \{y(j)\} $
of length $ n $ so that $ \sum_{j=1}^n y^2(j) = 1; $ applying the Lagrange 
factor method for calculation we find that the maximum is attained, in 
particular, on non - negative sequences, some component of which equals zero, 
while the positive components are equal. \par
 For instance if $ \ \forall |\lambda| \ge 1 \ \ \varphi_m(\lambda) = |\lambda|^q $
for some $  q = const > 1, \ $ then 
$ \ \chi_m(\lambda) \sim C |\lambda|^{\max(2,q)}, \ |\lambda| 
\to \infty. $  \par 
{\bf Lemma 2.} {\it Let the variables} $ \{\xi(i,m)\} $ {\it  be centered, independent, 
and assume that the condition (2.7) is satisfied, then }
$$
\sup_{b \in B} T \left(\sum_i b(i) \xi(i,m), x \right) \le \exp \left(-\chi_m^*(x)\right), \eqno(2.9)
$$
{\it where} 
$$
\chi^*(x) = \sup_{\lambda} (\lambda x - \chi(\lambda)) \ -
$$
{\it is the so - called Young - Fenchel transform.}  \par
{\bf Proof}. Taking into account the independence and  assuming $ \lambda > 0 $ we obtain:
$$
{\bf E} \exp (\lambda \sum_{i=1}^n  b(i) \xi(i,m)) = \prod_{i=1}^n {\bf E} \exp(\lambda b(i) 
\xi(i,m)) \le 
$$
$$
\le \exp \left( \sum_i \varphi_m(\lambda b(i)) \right) \le \exp 
\chi_m(\lambda)
$$
 by virtue of definition of the function $ \chi_m(\cdot). $ The 
proposition  of  lemma 2 it follows now from Chebyshev inequality. See in detail 
 [11, p.24.]\par
{\bf Corollary 1.} We obtain  syntesing propositions   of lemma 2 and lemma 1,  that
in our assumptions 
$$
\sup_{b \in B} T \left(\sum_i b(i) \xi(i,m), \ x \right) \le \min \{W[T_m](x), 
\exp \left(-\chi^*_m(x) \right) \} \stackrel{def}{=}
$$
$$
\stackrel{def}{=}  \overline{W}[T_m](x).
$$ 
 For instance suppose that $ T(\xi(i,1),x) \le \exp \left(-x^q \right), \ q > 0, $
and  recall that the function $ N(q) $ at $ q \in (0,1] $ is  equal to 
$ 2q/(2+q) $ and at $ q \in (1,\infty] \ \Rightarrow  \ N(q) = 
\min(q,2). $  Then we get at $ x \ge 1 $
$$
\sup_{b \in B} T(\sum_i b(i) \ \xi(i,1), x)  \le \exp \left( - C_7(q)
x^{N(q)} \right). 
$$
 It is  proved in [11, p. 50] that the exponent $ \min(q,2) $ in the case $ q > 1 $ is 
non - improvable. \par
 More generally assume that $ \exists q > 0, r \in (-\infty,\infty) \ \Rightarrow \forall x > 0 $
$$
T(\xi(i,1),x) \le \exp \left(-C x^q (\log(F(q,r) +x))^r \right).
$$
 Introduce the following vector - function $ N(q,r) = (N(1;q,r),N(2;q,r)): $ \\ 
a) at $ q \in(0,1) $ or $ \ q=1, r < 0 \ \Rightarrow $
$$
  \Rightarrow N(1;q,r)=2q/(q+2), \ \  N(2;q,r) = 2r/(q+2);
$$
b) at $  q=1, r \ge 0 $ or $ q \in (1,2), $ or $ q=2,\  r < 0 \ \Rightarrow $   
$$
 \Rightarrow N(1;q,r)=q, \ \ N(2;q,r)=r;
$$

c) at $ q=2, r \ge 0 $ or $ q > 2 \ \Rightarrow $
$$
 \Rightarrow N(1;q,r)=2, \ \  N(2;q,r)=0. 
$$
Proposition:  $ \ \sup_{b \in B} T(\sum_i b(i) \ \xi(i,1),x) \le $ 
$$
 \exp \left(-C_8(q,r) x^{N(1;q,r)} \ (\log(F(N(1;q,r),N(2;q,r)) + x))^{N(2;q,r)} \right).
$$
{\bf Theorem 5.} {\it Assume in addition that the r.v.} $ \{\xi(i,m) \} $ {\it are independent.
Let us define the sequence of tail - functions by the following initial condition and recursion:}
$$
T^{(1)}(x) = \overline{W}[T_d](x),
$$
$$
T^{(m+1)}(x) = W \left[T_{m+1} \vee T^{(m)}  \right](x), \ m = 1,2,\ldots, d-1.
$$
 {\it Proposition: }
$$
\sup_{b \in B} T(Q_d, x) \le T^{(d)}(x). \eqno(2.10).
$$ 
{\bf Proof}. For the sum 
$$
Q(d,n) = \sum_{I \in I(d,n)} b(I) \prod_{m=1}^d \xi(i_m,m),
$$
i.e. without diagonal members, in the case $ d = 1 $ our result is the content 
of lemma 2. The recursion is proved likewisely the proof of theorem 4. \par
{\bf Proof} of  theorem 2. We must  prove that the non - diagonal members in 
expansion for $ R_d $ have "little" tails:
$$
\sup_{b \in B} T(R_d - Q_d,x)\le \exp \left(-C(q) x^{U+\varepsilon} \right), \ x \ge 1,\eqno(2.11)
$$
for some $ \varepsilon = const \ge 0 $ as long as $ T(\xi+\eta,x) \le T(\xi,x/2) +
T(\eta,x/2). $  The unique non - trivial case is if the family $ \{\xi(i,m)\},
 \ m \le d-2 $ does not depend on the $ \{ \xi(k,d-1) \} $ and $ \xi(i,d-1) =
\xi(i,d). $ Let us denote
$$
A  = \sum_{ I \in J(d,n)} b(I) \xi(i_1,1)\xi(i_2,2)\ldots, \xi(i_{d-2},d-2) 
[\xi^2(i_d,d) - m(2,i_{d},d)].
$$
We receive: 
$$
T(\xi^2(i,d) - m(2,i,d),x) \le \exp \left(- C(d)x^{q/2} \right), \ x \ge 1.
$$
 Thus  
$$
T(A, x) \le \exp \left(- C x^{N(q(1),q(2),\ldots,q(d-2), q/2)} \right), \ 
q = q(d).
$$
As long as
$$
1/N(q(1),q(2), \ldots,q(d-2), q/2) = (d-2)/2 + \sum_{j=1}^{d-2}
[1/q(j)] + 1/L(q/2),
$$
$$
1/N(q(1),q(2),\ldots,q(d-2),q,q) = (d-1)/2 + \sum_{j=1}^{d-2} 
[1/q(j)] + 1/q + 1/L(q),
$$
It is enough  to prove  the  following inequality: $ \forall q > 0 \ \Rightarrow $ 
$$
1/2 + 1/q + 1/L(q) - 1/L(q/2) \ge 0.\eqno(2.12)
$$
 After the consideration of all cases $ q \in (0,1], \ q \in (1,2], \ q \in 
(2,4] $ and $ q > 4 $ we can see  that the inequality (2.12) is valid. \par

 {\it Low bounds for tails in the martingale cases (theorem 1).}  Assume  
$ d = 2; $  the general case provides analogously. \par
 Let us introduce two martingale - differences:
$$
\xi(i) = \tau \nu(i), \ \eta(j) = \rho \theta(j),
$$
where all variables $ \tau, \{\nu(i) \}, \rho, \{\theta(j)\} $ are independent and
$ \{\nu(i) \}, \{\theta(j)\} $ are Rademacher series: 
$$
{\bf P}(\theta(j) = 1) = {\bf P}(\theta(j) = - 1) = {\bf P}(\nu(i) = 1) = 
{\bf P}(\nu(i) = - 1) = 1/2
$$
and $ \tau \ge 0, \ \rho \ge 0, $
$$
T(\tau,x) = \exp \left(-x^{q(1)} \right), \ T(\rho,x) = \exp 
\left(-x^{q(2)} \right),
$$
 and let us introduce at $ n \ge 2 $  the probability 
$$
{\bf P}_n(x) = {\bf P} \left( \sum_{1 \le i < j \le n} \tau \ \rho \
\nu(i) \ \theta(j)/ \sqrt{n(n-1)} > x \right), \ x > 1.
$$
We see that for  $ \ y,z \ge 1, \ x/(yz) > 1: $
$$
{\bf P}_n(x) \ge {\bf P} \left(  \sum \nu(i) \ \theta(j)/\sqrt{n(n-1)} > x/(yz) \right)
\exp \left(-y^{q(1)} \right) \exp \left(-z^{-q(2)} \right).
$$
 As long as (see [24])
$$
\lim_{n \to \infty} {\bf P} \left( \sum_{1 \le i< j \le n} \nu(i) 
\theta(j)/\sqrt{n(n-1)} > x/(yz) \right) = {\bf P} \left(I(g) > x/(yz) \right),
$$
 where $ I(g) $ is a two - dimensional stochastic integral (1.11) with some non - trivial function 
$ g(\cdot) \in L_2(R^2).$  Since 

$$
{\bf P} \left(I(g) > x/(yz) \right) \ge \exp \left(-C(g) x/(yz) \right),
$$
we obtain:
$$
\overline{\lim}_{n \to \infty} {\bf P}_n(x) \ge \exp \left(- C(g) x/(yz) - y^{q(1)} -
z^{q(2)} \right). \eqno(2.14)
$$
 Since
$$
 \sup_n {\bf P}_n(x) \ge \overline{\lim}_{n \to \infty} {\bf P}_n(x),
$$
 after the maximization the right - side of inequality (2.14) over $ y,z, $ we deduce
the left inequality in the theorem 1.\par 
 In order to prove the low bounds for tails distributions in the independent 
case  (theorem 3), we must consider  two examples. 
Taking  $ Q(d,n) = \prod_{i=1}^d \xi(i),$ where $ \xi(i)  $ are independent, 
symmetrical and such that 
$$
T(\xi(i), x) = \exp \left(-x^q \right), \ x \ge 0,
$$
we deduce after some calculations:
$$
T(Q(d,n), x)  = {\bf P}(\prod_{i=1}^d \xi(i) > x) \ge \exp 
\left(- C(d,q) x^{q/d} \right).
$$
 The sequence $ Q(n,d) $ converge in distribution at $ n \to \infty $
 by appropriate choice of  {\it sequences } of coefficients $ b(I) = b(I;n) $  
 to the $ d \ - $ multiple stochastic integral 
over the Gaussian orthogonal stochastic measure $ Z(\cdot) $ \ (1.11):
$$
Q(d,n)/\sqrt{ {\bf D} \ Q(d,n)} \stackrel{d}{\to} I(h) = \int_{R^d} h(y) Z(dy),
$$
with tail behavior $ T(I(h), x) \ge \exp \left(-C x^{2/d} \right). $ Consequently, 
$$
\sup_{b \in B} \sup_{\{\xi(i,m)\}: ||\xi(i,m)||_{q(m)} \le 1 } T(Q_d,x) \ge \exp 
\left(-C x^{\min(2,q)/d} \right). 
$$

\vspace{3mm}

{\bf 3. Moment estimations.} \ We shall derive here the moment estimations for $ Q_d $ 
in martingale and independent cases in the 
terms of {\it unconditional } moments, more exactly, $ L_p $ norms: $ |\eta|_p = {\bf E}^{1/p}
|\eta|^p $  of summands:
$$
\mu_m(p) \stackrel{def}{=} \sup_i {\bf E}^{1/p} |\xi(i,m)|^p = \sup_i |\xi(i,m)|_p, \ p \ge 2.
$$
 We define the sequence $ \gamma(d), d=1,2,\ldots $ by the following  
initial condition: $ \gamma(1) = \sqrt{2} $ and recurrent equation:
$$
\gamma(d+1)= \sqrt{2} (1+1/d)^d \gamma(d). 
$$
 In particular, $ \gamma(2)=4, \gamma(3) = 9 \sqrt{2}.$
Since $ (1+1/d)^d < e, $ we conclude: $ \forall d \ge 3 \ \Rightarrow \gamma(d) \le 9 
\exp\left((d-3)+(\log 2)(d-2)/2 \right).$ \par
{\bf Theorem 6.} {\it If} $ \ \forall m=1,2,\ldots,d \ \ \mu_m(d\cdot p)< \infty, $ {\it then}
$$
\sup_{b \in B}|Q_d|_p \le \gamma(d) \ p^d \ \prod_{m=1}^d \mu_m(d\cdot p).\eqno(3.1)
$$
{\bf Proof.} We start our consideration from the case $ d=1.$ The beginning is like one in 
[10],  but further instead of  convexity we intend to employ  $  H\ddot{o}lder \ \ inequality. $ \par 
 Note that it can be assumed that  $ p > 2 $ (the case $p=2$ is trivial) and 
 $ \forall i \le n \ \Rightarrow \ b(i) \ne 0.$ Further, the sequence $b(i) \ \xi(i), \
\xi(i) \stackrel{def}{=} \xi(i,1) $ is also the sequence of martingale differences relative to the 
initial filtration. Using the Burkholder inequality ([20], p. 78 - 81; [31])  we obtain
$$
|\sum b(i) \ \xi(i)|_p^p \le (p\sqrt{2})^p \ {\bf E} \left[\sum b^2(i) \ \xi^2(i) 
\right]^{p/2}.\eqno(3.2)
$$
 Let $ a(i) $ be some positive non - random sequence, the choice of which we shall  
clarify below. By virtue of $ H\ddot{o}lder \ \ inequality,$ in which we substitute 
$$
\beta = p/2, \ \  \alpha = \beta/(\beta-1) = p/(p-2),
$$
we see get
$$
\sum b^2(i) \xi^2(i) = \sum a(i) b^2(i) \cdot [\xi^2(i)/a(i)] \le 
$$
$$ 
 \left[\sum (b^2(i) a(i))^{\alpha}  \right]^{1/\alpha}  \cdot \left[\sum (\xi^2(i)/a(i))^{\beta} 
\right]^{1/\beta}.
$$
 Let us choose $ a(i) $ by the formula 
$$
a(i) = |b(i)|^{-4/p},
$$
then we receive 
$$
\left[\sum(b^2(i) a(i))^{\alpha}\right]^{p/(2\alpha)} = \left[ \sqrt{\sum b^2(i)} \right]^{p-2};
$$   

$$
{\bf E}\left[ \sum(\xi^2(i)/a(i))^{\beta} \right]^{p/2\beta} \le {\bf E}\sum |\xi(i)|^p b^2(i) \le 
\mu_1(p) \sum b^2(i).
$$
 Substituting the last inequality into (3.2) we obtain the proposition of theorem 6. \par

{\bf Corollary 2.} Let us denote 
$$
K_M(1,p)=\sup |\sum_{b \in B} b(i)\xi(i)|_p/\mu(p),
$$
where the upper bound is calculated over all sequences of centered martingale differences 
$ \{\xi(i)\} $ with finite absolute moments of order $p$. It follows from [3.1] that $ K_M(1,p) \le 
p\sqrt{2}. $ But it is proved in [13], [22] etc. that  if $ \{\xi(i)\} $ are independent 
symmetrical identically distributed, then the fraction 
in the right-hand part can has an estimate from below of the form $ 0.87 p/\log p. $ Thus
$$
0.87 p/\log p \le K_M(1,p) \le \sqrt{2} \ p.
$$
 Therefore  our estimation can'not be essentially improved.\par
{\bf Remark 1.} It the {\it other terms} of conditional expectations  $ {\bf E} \xi^2(i,m)/F(i-1), \
 {\bf E} \max_{i \le n} |\xi(i,m)|^p  $ the upper bound for the constant which is like to our 
$ K_M(1,p) $  is obtained in [5] then growth  by $ p \to \infty $ like $ p/\log p.$  
Our result can not be obtained from [5].\par
 The finall estimate for $ K_M(1,p),$ i.e. the  exact grows for $ K_M(1,p) $ by $ p \to \infty $
is now unknown.\par 
 We shall prove the general case $ d > 1 $ by induction over  $d$. We get for the corresponding 
martingale differences $ \zeta(n)=\tau(n,d)\eta(n,d) $ of a martingale $ (Q(n),F(n)) $   
on the basis  $  H\ddot{o}lder \ \ inequality $ and induction statement:
$$
|\zeta(n)|_p \le |\eta(n,d)|_{pd} \ |\tau(n,d)|_{pd/(d-1)} = \mu_d(pd) \cdot|Q(d-1)|_{pd/(d-1)} \le 
$$
$$
 \mu_d(pd) \ (pd/(d-1))^{d-1} \gamma(d-1) \sqrt{ \sum_{I \in J(d)} b^2(i) }  \prod_{m=1}^{d-1}
\mu_m(p(d-1)d/(d-1))=
$$

$$
 \gamma(d-1)[pd/(d-1)]^{d-1}\prod_{m=1}^d \mu_m(pd) \ \sqrt{ \sum_{I \in J(d)} b^2(I) }.
$$
  We receive  from the one - dimensional case:
$$
|Q(d,n)|_p \le p \sqrt{2} \ \gamma(d-1) \ [pd/(d-1)]^{d-1} \prod_{m=1}^d \mu(pd) \times
$$
$$
 \sqrt{ \sum_{I \in I(d,n)}  b^2(I)} = \gamma(d) \ p^d \ \prod_{m=1}^d \mu_m(pd).
$$
 (Recall that  $ b \in B \ \Leftrightarrow \ \sum_{I \in I(d,n) } b^2(I) = 1.$).\par 
 {\bf Corollary 3.} Let us denote
$$
K_M(d,p) = \sup_n \sup_{ \{\xi(i,m)\}: |\xi(i,m)|_{d \cdot p} < \infty} |Q_d|_p /
\prod_{m=1}^d \mu_m(d \cdot p),
$$
where $ \sup $ is calculated over all families of sequences of centered martingale differences
$ \{ \xi(i,m) \} $ with condition $ |\xi(i,m)|_{d \cdot p} < \infty. $ It follows  
 from theorem 6 that $ K_M(d,p) \le \gamma(d) \ p^d. $ We shall prove now 
{\it the low bounds for } $ K_M(d,p): $ at $ d \ge 2 $
$$
K_M(d,p) \ge C^d \ p^{d/2},
$$
here $ C  $ is an absolute constant. Let us choose the family of independent Rademacher series:
$$
{\bf P}(\xi(i,m) = \pm 1) = 1/2,
$$
then $ \mu_m(d \cdot p) = 1. $ Introduce at  $ n \ge d $ the random variables
$$
\overline{Q}_d(n) = \overline{Q}_d = \sum_{I \subset I(d,n)} \prod_{m=1}^d 
\xi(i_m,m) /[n(n-1)\ldots (n-d+1)].
$$
 We conclude relaying on the main result of the article [24]: at $ n \to \infty \ \Rightarrow 
\ \overline{Q}_d \stackrel{d}{\to} I(h) $ in the sense of distribution convergence,
where $ I(h) $ is  the  multiply stochastic integral (1.11) with tail behavior 
$$
T(I(h),x) \ge \exp \left( - C_1 x^{2/d} \right), \ x \ge C_2.
$$
 It is easy to verify that the moment convergence in [24] is  true. Hense 
$$
K_M(d,p) \ge \overline{\lim}_{n \to \infty} |\overline{Q}_d|_p \ge C_3^d \ p^{d/2}.
$$

 Now we shall consider again "independent case", under condition: for some $ p \ge 2 $
$$
\mu_m(p) \stackrel{def}{=} \sup_{i \le n} |\xi(i,m)|_p < \infty, \ m=1,2,\ldots,d.
$$

{\bf Theorem 7.} {\it If all the variables} $ \{\xi(i,m)\} $ {\it are independent, centered 
and so that  for some} $ p \ge 2 \ \Rightarrow \mu_m(p) < \infty \ \ \forall m = 1,2,\ldots,d, $ 
{\it then}
$$
\sup_{b \in B} |Q_d|_p \le 2^{d/2} p^d \left[\prod_{m=1}^d \mu_m(p)\right]/\log p. \eqno(3.3)
$$
{\bf Proof} is the same as in the theorem 6. In  the case $ d=1 $ the proposition (3.3) is 
provided in [13]; we can rewrite this result on the form 
$$
|\sum b(i) \xi(i,1)|_p \le \sqrt{2} \ \sqrt{ \sum_{i=1}^n b^2(i)} \ \mu_1(p)/\log p, \ \ p \ge 2.
$$
  Since the variables $ \{\xi(i,m)\} $ are independent  we see that for $ \zeta(n) = 
\eta(n,d) \ \tau(n,d) $ the following inequality holds: 
$$
|\zeta(n)|_p = |\eta(n,d)|_p \ |\tau(n,d)|_p \le \mu_d(p) |Q(d-1,n)|_p.
$$
 It follows from induction statement that: 
$$
|Q(d,n)|_p \le 2^{d/2} \ p \ \mu_d(p) \ p^{d-1} \ \prod_{m=1}^{d-1} \sqrt{\sum_{I \in I(d,n)} b^2(I)
}/\log p =
$$
$$
 2^{d/2} \ p^d \ \prod_{m=1}^d \mu_m(p) \ \sqrt{\sum_{I \in I(d,n)} b^2(I)}/\log p.
$$
 Note that our result (3.3) improves the same estimations in article [35].\par
 {\bf Corollary 4.} Let us denote {\it in the independent case} 
$$
K_I(p,d) = \sup_n \sup_{b \in B}\sup_{\{\xi(i,s): |\xi(i,m)|_p < \infty \}} |Q(d,n)|_p /
\prod_{m=1}^d \mu_m(p).
$$
 Proposition: for some absolute constant $ C_4 \in (0,\infty) $ 
$$
C^d_4 \ p^d \log^{-d} p \le  K_I(p,d) \le 2^{d/2} \ p^d /\log p.   
$$
 {\bf Proof} the low bounds: the moment estimations are derived in [21], [29] for 
symmetrical polynomials on independent identical symmetrically distributed  variables:
$$
|Q(d,n)|_p /\left(\sqrt{ {\bf D}Q(d,n)} \ \mu^d(p) \right) \ge C_4^d p^d /(\log^d p). 
$$
{\it Our hypothesis:} $ p \to \infty \ \Rightarrow K_I(p,d) \asymp C(d) \ p^d/(\log p)^d.$\\
{\bf Corollary 5.} We want  make the comparision between exponential and moment estimations
in the independent case. Let us introduce, more exactly,
 for any Orlicz space of random variables $ W $ with norm $ ||\eta||W $  
on the our probability space  the uniform tail of distribution
$$
{\bf P}_W(x) = \sup_{\{||\xi(i,1)||W \le 1\}} \sup_{b \in B} T(Q_1,x).
$$
 There are two methods of estimation of $ {\bf P}_W(x): $ "exponential estimations"  and "moment 
estimations". Namely, if a norm $ ||\cdot||W $  in a $W$  space is equivalent to the norm of a view 
$$
|||\eta|||W = \sup_{p \ge 2}|\eta|_p /\psi(p),
$$
where $ \psi(\cdot)  $ is a monotonically increasing continuous 
positive function, $ \psi(+\infty) = \infty, \ $ 
for example  norms in the spaces $ L_p, G(q), G(q,r), \Psi(\beta), $ then by means of theorem 7 
we can estimate  all moments of $ Q_1 $ and further by means of Chebyshev inequality
the tail of distribution $ Q(d,n). $ For the spaces $ L_p, p \ge 2 $ the 
moment method gives better result; in the case $ W = \Psi(\beta) $ both  methods give the same 
result:
$$
{\bf P}_W(x) \le \exp \left(-C(\beta) (\log(1+x))^{1+1/\beta} \right); 
$$
for the spaces $ G(q,r) $ the  the exponential method is better.\par
 Let us consider the following example. Assume that for some  $ r = const, q = const > 0 $
$$
\sup_i \sup_m T(\xi(i,m),x) \le \exp \left(-x^q (\log(F(q,r)+x)^r) \right), \ x>0.
$$
 Denote $V(q,d)= 2q/(d(q+2)) $ and define the sequence $ r(d) $ by the  recurrent equation:
$$
r(d+1)= [r(d)q + V(d,q)r]/[V(d,q)q+2V(d,q)+2q]
$$
with initial condition $ r(1)=2r/(q+2). $ Then  in the  martingale case 

$$
\sup_{b \in B} T(Q_d,x) \le \exp \left(-C(d,q) x^{V(d,q)} (\log(F(V(d,q),r(d)) + x)^{r(d)}\right).
$$
If $ \sup_i \sup_m T(\xi(i,m) \le \exp \left(-(\log(1+x))^{1+1/\beta} \right), \ \beta >0,$
then 
$$
\sup_{b \in B} T(Q_d, x) \le \exp \left( - C_1(d,\beta) (\log(1+x))^{1+1/\beta} \right).
$$
 {\bf Remark 2}. We can derive similarly using Doob's inequalities 
the  moment and exponential estimates for the distribution  of the variables
$$
\max_{k=1,2,\ldots,n}  Q(d,k), \ \max_{k=1,2,\ldots,n}  |Q(d,k)|,
$$
 {\bf Remark 3.} It is easy to receive the generalization of our inequalities in the so - called 
martingale fields $ \xi(i,m) = \xi(\vec{i},m); $ see definitions and some preliminary results  
 in [10].\par 

\vspace{3mm}
{\bf 4. Application to the theory of $ U \ - $ statistics.} \
 In this section we shall apply  our estimations ( (2.6), (3.3) etc.) in the theory 
of $ U \ - $ 
statistics. For more detail about  using  of martingale technique based on the 
Hoeffding martingale representation for $ U \ - $ statistic see in [37], [16]. 
 Our results improve also somewhat the estimations in [16], [21] etc.\par 
 Let $ \{\xi(i)\}, \ i = 1,2,\ldots,n $ be independent identically distributed random variables with 
values in the fixed measurable space $ \{ X,S \}, \ \Phi(x_1,x_2,x_3,\ldots,x_d), \ d < n  $ be 
a symmetrical measurable non - trivial numerical function (kernel)  of $ d $ variables:
 $ \Phi: X^d \to R^1, $

$$
  U(n) =  U(n; \Phi,d) = 
 \sum_{I \in I(d,n)} \Phi(\xi(i(1),\xi(i(2)), \ldots, \xi(i(d))) / {n \choose d}
$$
be a so - called $ U \ - $ statistic. Denote $ Dim \ \Phi = d, $
$$
\Phi = \Phi(\xi(1),\xi(2), \ldots,\xi(d)), \ r = rank \ U \in [1,2,\ldots,d-1];  
$$
$$
T(\{\Phi,d\},x) \stackrel{def}{=} \sup_{n > d} T((U(n)/\sqrt {{\bf D} U(n)}, \ x).
$$

 Assume that $ {\bf E}\Phi = 0, {\bf D} \Phi \in  (0,\infty) $ and  all the   moments  
$ \Phi $ which is written below there  exist; otherwise the results are trivial. Recall
 here for readers convenience  the so - called 
{\it martingale representation } for the centered $ U \ - $ statistic  (see [16], p. 26; [37]): 
$$
U(n) = \sum_{k=r}^d {k \choose d} U(n;k),  \ \ U(n;k) = \sum_{I \in I(n,k)} 
g_k(\vec{\xi}(I))/ {n \choose k},
$$
where $ \vec{\xi}(I) = \{\xi(i(1)), \xi(i(2)), \ldots, \xi(i(d)) \}, $ 
$ \mu(A) = {\bf P}(\xi(i) \in A), \ \int_X f(y) \delta_x(dy) = f(x), \  g_k(x_1,x_2,\ldots,x_d) = 
\Phi_k = g_k[\Phi](x_1,x_2, \ldots, x_d) = $
$$
= \int_{X^d} \Phi(y_1,y_2, \ldots,y_d) \prod_{l=1}^k \left( \delta_{x_l}(dy_l) - \mu(dy_l) \right)  \ 
\prod_{l=k+1}^d \mu(dy_l).
$$
 It is well known ([16], [37]) that the sequence $ S(k) = S^{(n)}(k) = C(k,n) \ U(n,k), \
 k \le n  $ relatively {\it some filtration } $ F(k) = F^{(n)}(k) $ is a martingale:
$$
{\bf E} S(l)/F(k) = S(k), \ k \in [1,l],
$$
 and that by $ n \to \infty \ \Rightarrow \ {\bf D} U(n) \asymp n^{-r}. $ \par  
{\bf Theorem 8.}  
$$
|U(n)/\sqrt{{ \bf D}U(n)}|_p \le  C^{d} \ p^{d} \ |\Phi|_p \ /\log p. \eqno(4.1)
$$
{\bf Proof.} The case $ d = 1 $ was consider in the section 3; now we shall use the method of   
induction over $ d;$ for simplicity we shall investigate only the case $ d = 2.$ \par
 Assume for beginning that the (non - trivial) kernel $ \Phi $ is  degenerate:
$$
{\bf D} \Phi(\xi(1), \xi(2))/\xi(1) = 0 \ (mod \ {\bf P }), \ (\Leftrightarrow r \ge 2).
$$

 The sequence 
$$ 
 {n \choose 2} U(n)  =  \sum_{1 \le i<j \le n} \Phi(\xi(i),\xi(j))
$$ 
is a martingale relative to  {\it some }  $ \sigma - $ flow (filtration) 
 with correspondence martingale - differences    
$$
\zeta(n) = \sum_{i=1}^{n-1} \Phi(\xi(i), \xi(n)),
$$
  Fixing the value of $ \xi(n), $ and  using the induction statement,  and denoting 
$ \mu(A) = {\bf P}(\xi(i) \in A), $ at $ p,n \ge 2, $ we obtain: 
$$
|\zeta(n)|^p_p = {\bf E} |\zeta(n)|^p = \int_X \mu(dx) {\bf E} \left|\sum_{i=1}^{n-1} 
\Phi(\xi(i),x)\right|^p \le
$$

$$
 \int_X \mu(dx) \ (n-1)^{p/2} \ (p\sqrt{2}/\log p)^p \ {\bf E}|\Phi(\xi(1),x)|^p =  
$$

$$
 (n-1)^{p/2} \ (p \sqrt{2} /\log p)^p \ |\Phi|_p^p. 
$$
Hence 
$$
|\zeta(n)|_p \le \sqrt{2n-2} \ (p/\log p) \ |\Phi|_p.
$$ 
 We can  prove our proposition in  the case degenerate kernel
 substituting into the inequality (1.6) the least estimation. \par
 Now we shall consider the general (and non - trivial)  case $ r=1.$ We denote 
$ {\bf E} \Phi(\xi(1), \xi(2))/(\xi(2)=x) = g(x),$ 
$$
\Phi^0(x,y) = \Phi(x,y) - g(x) - g(y) = g_1[\Phi](x,y).
$$
 Then $  {\bf E} g(\xi(i)) = 0 $ and $ \Phi^0 $ is a degenerate kernel. It follows 
from Iensen inequality for 
conditional expectation  that $ |g(\xi(i))|_p \le |\Phi|_p,$   and, hence 
$ |\Phi^0|_p \le C(d) |\Phi|_p.$ \par

 Let us write the  Hoeffding decomposition for $ U \ - $ statistic: $\sqrt{n} \ U(n) = $ 
$$
 2/[ \sqrt{n} \ (n-1)]  \sum_{1 \le i < j \le n} \Phi^0(\xi(i),\xi(j)) +
 [4 \sqrt{n})/(n-1)] \cdot  \left( \sum_{i=1}^{n-1} g(\xi(i)) \right) -
$$
$  2 g(\xi(n)) \ \sqrt{n} /(n-1). $
 Using triangle inequality and our proposition for degenerate statistics, we obtain:
$$
| \sqrt{n} \ U(n)|_p \le C_1 \ n^{-1/2} \ p^2 \ |\Phi|_p /\log p + C_2 \ p \ |g(\xi(i))|_p /\log p \le 
$$
$$
 C_3 \ p^2 \ |\Phi|_p \ /\log p; 
$$
here $ C_{1,2,3} = C_{1,2,3}(d). $ \par 

 We shall derive now the exponential bounds for tail of distribution $ T( {\Phi,d},x)). $ 

{\bf Theorem 9.} {\it Assume that for some} $ K \in (0,\infty), q > 0, r \in R^1 \ \Rightarrow  $
$$
T(\Phi,x) \le \exp \left( -(x/K)^q \  (\log (1+x/K))^{-qr} \   \right), \ x \ge 0. \eqno(4.2)
$$
{\it Then}  $ T(\{\Phi,d \},x) \le $
$$
 \exp \left(  - C(d,q,r) (x/K)^{q/(qd+1)} 
(\log(1+x/K))^{-((r-1)q)/(qd+1)} \right). \eqno(4.3)
$$
{\bf Proof}. We can assume  $ K = 1. $  It follows from  condition (4.2) that:
$$
|\Phi|_p \le C_1 \ p^{1/q} \ \log^r p, \ \ p \ge 2.
$$
 We conclude using theorem 8: 
$$
| U(n)/\sqrt{ {\bf D} U(n) }|_p \le C_2 \ p^{d+1/q} \ \log^{r-1} p. 
$$
 We deduce  (4.3) returning to the tail of probability \par
 Now we shall receive the refined  but more cumbersome 
exponential bounds for tail distribution $ T(\{\Phi,d \},x). $ We shall  
deduce the recurrent equation (on the dimension $ d = Dim \ \Phi $) for one in the 
spirit of the section 2. Let us denote 
$$ 
t(d,k,r) = 1/ \left[(d-r+1) \ {d \choose k}  \right]  
$$ 
and for almost all $ (mod \ \mu) $ values $ \ z \in X: g_{k,(z)}[\Phi_k] = $
$$
 \Phi_{k,(z)} = \Phi_{k,(z)}(x_1,x_2,\ldots,x_{d-1}) \stackrel{def}{=} 
\Phi_k(x_1,x_2,\ldots,x_{d-1},z).
$$ 
 Note  that 
$$
T(\{\Phi,d \},x) \le \sum_{k=r}^d T(\{ g_k[\Phi],d \}, \ t(d,k,r) \ x ). \eqno(4.4)
$$
  Consequently, it is sufficient to estimate this distribution only for degenerate kernels
$ \Phi_k(\cdot), $ i.e.  if $ r = rank \ U \ge 2, $ as long as the arbitrary kernel 
$ \Phi(\cdot) $ may be represented as a linear combination of degenerate kernels $ \Phi_k $
 ([16], p. 26). \\
{\bf Theorem 10.} {\it The tails} $ T(\{\Phi_k,d \},x) = T(\{g_k[\Phi],d \}, \ x) $ 
{\it may be estimated as} 
$$
T(\{g_k[\Phi],d \}, x) \le L(\{g_k[\Phi],d \},x ), \eqno(4.5)
$$
{\it where the functions} $ L(\{g_k[\Phi],d \},x) $ {\it satisfy the following system of 
a recurrent equations:}
$$
L(\{g_k[\Phi],d \},x) = W \left[ \int_X \mu(dz) \ L(\{g_{k,(z)}[\Phi],d-1 \},x) \right], \ d \ge 2, 
\eqno(4.6)
$$

{\it with initial condition} 
$$
L( \{\Phi,1 \},x) = \overline{W}[T(\Phi)](x). \eqno(4.7)
$$
{\bf Proof} (briefly; in  the case $ d = 2 $). Define the sequence of functions 
$ L(\{g_k[\Phi],d \},x) $ by means of the equations 
(4.6) and (4.7). Then the estimation (4.5) by $ d=1 $  it follows from the definition of operator that 
$ \overline{W}(\cdot), $ as long as the variables $  \{\Phi(\xi(i)) \} $ are independent. \par
 Further, since the kernel $ g_k[\Phi] $ is degenerate, we deduce that the sequence 
$$
\Delta(n) = \sum_{I \in J(d,n)}  \Phi_k(\vec{\xi}(I),z)
$$
for almost all values $ z, \ z \in X $ 
with respect to some filtration is a  martingale with correspondence martingale - 
differences  $ \zeta(n). $ 
Using the proof of theorem  4  and omitting some calculations we obtain 
  the inequality (4.5).\par

\vspace{3mm}
{\bf 5. Applications to the stochastic integration.} \ 
A). Let $ (Z(t),F(t)), \ t \ge 0 $ be a {\it left } - continuos in the $ L_2(\Omega,{\bf P})$ sense 
centered  square integrable martingale relatively the flow of $ \sigma \ - $ fields (filtration) 
$ F(t): F(0+0) = F(0) = \{ \emptyset, \Omega \}, \ F(t-0) = F(t), t >0,  Z(0) = 0. $ 
Denote for $ 0 \le a < c \le \infty \ \ \nu([a,c)) = 
{\bf D}(Z(c)-Z(a)); $ then $ \nu(\cdot) $ may be continued to the measure (may be unbounded) 
on the Borel subsets on half - line $ [0,\infty).$  \par
 Let $ b(t), \ t \ge 0 $ be a non - random measurable function belonging to the 
space $ L_2(R^1_+,\nu):$
$$
 ||b(\cdot)||^2(L_2(\nu)) \stackrel{def}{=} \int_{[0,\infty)} b^2(t) \ \nu(dt)  < \infty.
$$
 We  define for some $ q > 0 $ 
$$
 ||Z||(Lip(q),\nu) = \sup_{ \nu([a,c)) \in (0, \infty)} ||Z(c)-Z(a)||_q /
\{\sqrt{\nu(c) - \nu(a)} \},
$$
and consider a stochastic integral 
$$
[b;Z] \stackrel{def}{=} \int_{[0,\infty) } b(t) d Z(t).
$$
{\bf Theorem 11}. {\it If for some} $ q > 0 \ \ ||Z||(Lip(q),\nu) < \infty, $ {\it then
for some } $ C(q) \in (0,\infty) $ 

$$
|| \ [b;Z] \ ||_{M(1,q)} \le C(q) ||b||(L_2(\nu)) \cdot ||Z||(Lip(q), \nu).  \eqno(5.1)
$$
{\bf Proof.} It is enough to prove  our  proposition (5.1) for a simple functions 
$ b(\cdot),$ i.e  of a view $ b(t) = $  
$$
 = \sum_{k=0}^K b(t(k)) \chi(t \in [t(k), t(k+1)), \ 0 \le t(0) < t(1) < \ldots < t(K+1) < \infty.
$$
 Here $ \{t(k)\} \ $ is any {\it non - random sequence}  such that  $ \nu(t(k+1)) - 
\nu(t(k)) > 0.$  For those functions we can write: $ [b;Z] =$
$$
\sum_{k=0}^K b(t(k)) \ \{Z(t(k+1)) - Z(t(k)) \} = \sum_{k=0}^K b(t(k)) \sqrt{\nu([t(k),t(k+1))} 
\times
$$

$$
 \{Z(t(k+1)) -Z(t(k)) \}/ \{ \sqrt{\nu([t(k), t(k+1))} \}. \eqno(5.2)
$$
 We deduce by virtue of theorem 1,
using the proposition of theorem 1 for the sequence $ b(i): = b(t(i)) 
\sqrt{\nu(t(i+1)) - \nu(t(i))} $ and 
choosing  a martingale differences $ \{\xi(i)\} $ as
$$ 
\xi(i) = \{Z(t(i+1) - Z(t(i)) ) \}/\{\sqrt{\nu([t(i),t(i+1))} \}:
$$
$$
||[b;Z]||_{M(1,q)} \le C(q) \sqrt{ \sum_{k=0}^K b^2(t(k)) \ \nu([t(k),t(k+1)) } \cdot \sup_i 
||\xi(i)||_q =
$$
$$
 C(q) ||b||(L_2(\nu)) \ ||Z||(Lip(q),\nu).
$$
 On the other hand,  We can write
 denoting $ (b;Z)_q = C(q) ||b||(L_2(\nu)) \times ||Z||(Lip(q),\nu) $ 
 for enough greatest values $ x > x_0 $ the inequality: 
$$
\exp \left(-[x/(C_1(b;Z)_q]^{M(1,q)} \right) \le \sup_{b(\cdot):||b||=1} \ 
\sup_{Z:||Z||(Lip(q),\nu)<\infty} T([b;Z],x) \le  
$$
$$
 \exp \left(-[x/(C_2(b;Z)_q)]^{M(1,q)} \right).
$$
 Low bounds  follow from the left inequality of Theorem 1 (1.5).\par
 It is easy to see that we can receive  analogous result for {\it multiply } $ d - $ 
dimensional stochastic integral of a kind $ [b; Z_1,Z_2,\ldots,Z_d] = $
$$
  \int \int \ldots \int_{0 \le t(1)<t(2)< \ldots t(d) <\infty}
b(t(1),t(2),\ldots,t(d)) \prod_{m=1}^d Z_m(dt(m)),
$$
where  $ b(t(1),t(2),\ldots,t(d)) $  are non - random measurable square integrable functions,
$ Z_m(t) $ are the square integrable centered left - continuous martingales with the 
correspondences measures $ \nu_m([a,c)) = {\bf D} \left(Z_m(c) - Z_m(a) \right) $ and values
$ ||Z_m||(Lip(q(m),\nu_m)): \ \exists C(\vec{q}) \in (0,\infty) \ \Rightarrow $

$$
|| \ [b; Z_1,Z_2,\ldots, Z_d] \ ||_{M(d,\vec{q})}  \le 
$$
$$
\le C(\vec{q}) ||b||L_2(R^d,\prod_{m=1}^d \nu_m) \ \prod_{m=1}^d ||Z(m)||Lip(q(m),\nu_m).
$$
 This result may be improved 
 in  the case when all martingales $ Z_m $ are independent and have independent 
increments  as in the theorems 2,3.  For Gaussian martingales $ Z_m(t) $ in ([15],  
p.119) there is more exactly result. But we  consider the other  problem. \par 
B). Let again $ (Z(t), F(t)), \ t \ge 0 $ be a square integrable centered {\it continuous } 
(with probability one) 
 martingale  with correspondence {\it quadratic variation }  $ <Z,Z>_t; \ 
<Z,Z>:=<Z,Z>_1. $ Let us consider the $ d \ - $ dimensional multiply stochastic 
integral $ [1;Z,(d)] = $
$$
 [1; Z,Z,\ldots, Z,(d)] = \int \int \ldots \int_{0 \le t(1)< t(2) \ldots < t(d) \le 1} 
\prod_{m=1}^d d Z(t(m)).
$$
 It is proved in the article  [39] the following generalization of the  classical 
Burkholder - Davis  - Gundy inequality (in our terms and notation):
$$
 | \ [1;Z,Z,...,Z] \ |_p \le  A(d,p) |<Z,Z>^{1/2}|_p,
$$
$$
 A(d,p) \stackrel{def}{=} (1+1/p)^d \ (dp)^{d/2}/d! \  \le C^d (p/d)^{d/2}, \eqno(5.3)
$$
where $ p \ge 2 $ and $ C $ is absolute constant.\par

{\bf Theorem 12}.  {\it Assume that  for some } $ q > 0 $
$$
||<Z,Z>||_q < \infty. \eqno(5.4)
$$
{\it Then }

$$
1. \ \ ||[1;Z,(d)]||_{M(d,q)} \le C(d,q) ||<Z,Z>^{1/2}||^d_q. \eqno(5.5) 
$$
{\it 2. If  (5.4)  holds for some} $ q > 2, $  {\it then taking  now} $ \beta = \beta(q) = 
(q+2)/(q-2) > 0 $ {\it we assert that the integral series}  $ 1+ \sum_d [1;Z,(d)] $ {\it for 
 Dolean exponent} $ E_D(Z)= \exp(Z(1) -0.5 <Z,Z>) $ {\it convergent in the} $ \Psi(C,\beta(q)), \ 
(\exists \ C \in (0, \infty)) $ {\it norm:}  
$$
|| E_D(Z) ||\Psi(\beta(q)) \le 1 + \sum_{d=1}^{\infty} ||[1;Z,(d)||\Psi(\beta(q)) < \infty. \eqno(5.6)
$$
{\bf Remark 4}. The case $ q=2 $ is considered in  [39]. \par
{\bf Proof}. We shall assume without loss of generality  $ ||<Z,Z>^{1/2}||_q \le 1/C_2, $ 
where $ C_2 $ is a "great" constant, so that 
$$
\forall p \ge 2 \ \Rightarrow |<Z,Z>^{1/2}|_p \le p^{1/q}.
$$
 We receive  using  (5.4) and Stirling formula 
$$
|[1; Z,(d)]|_p \le C^d \ p^{d/2} \ d^{-d/2} \ p^{d/q} \le C^d p^{1/M(d,q)}\ d^{d/q - d/2}. 
\eqno(5.7)
$$
Hence 
$$
||[1;Z,(d)]||_{M(d,q)} \le C^d d^{-d/2} < \infty, 
$$
 and  hence  we deduce (5.6). \par
 Further, assume  $ q > 2. $ We can write from (5.4) $ \forall p \ge 2:$
$$
| E_D(Z)|_p \le 1 + \sum_{d=1}^{\infty} |[1;Z,(d)]|_p \le C + \sum_{d=1}^{\infty} p^{d(0.5 + 1/q)} \ 
d^{-d(0.5 - 1/q)} \le
$$

$$
 \exp \left(C p^{(q+2)/(q-2)} \right) = \exp \left(C p^{\beta(q)} \right).
$$
 Now the estimation  (5.6) it follows from the definition of $ \Psi(\beta,C) \ - $ spaces.\par
 For instance assume that (5.4) holds at $ q = +\infty, $ i.e. \\
 $ vraimax_{\omega}|<Z,Z>| < \infty.$ Then $ \beta = 1, $  hence 
$$
T(E_D(Z),x) \le \exp \left(-C \log^2(1+x) \right), \ x > 0. \eqno(5.8)
$$

 In  the case when  $ Z(t) $ is the Wiener martingale, estimation (5.8) is exact as long as 
$ <Z,Z> = 1 $ and $ E_D(Z) = \exp(Z(1) - 0.5). $ \\ 

\vspace{2mm}

{\bf 6. Applications to the martingale summs.} \
 We shall concider in this section  the tail behavior of the centered  {\it martingale sums} (1.0)
$$
Q_d  = \sum_{I \in I(d,n)} b(I) \xi(I), 
$$ 
where $ n \le \infty, \ \sup_{i \le n} ||\xi(i,m)||_{q(m)} <\infty, $ or, equally, 
$  \exists K(m), q(m) \in (0,\infty] $ such that
$$
 \sup_{i \le n} T(\xi(i,m),x) \le \exp \left(- (x/K(m))^{q(m)} \right), \ x > 0,
\eqno(6.1)
$$ 
if the speed of convergence $ b(I) \to 0 $ by $ I \to Z_+^d $ is rapid. \par
 We obtain in the sections 2,3  the estimation of tail $ T(Q_d,x) $ under condition 
$ \sum b^2(I) < \infty. $  But if 
$$ 
\sum_{I \in I(d,n)} |b(I)| < \infty,  \eqno(6.2)
$$
it follows from triangular inequality for the $ G  \ - $ norms that
$$
|| Q_d ||_G \le C \sum_{I \in I(d,n) } |b(I)| \prod_{m=1}^d \sup_{i \le n} ||\xi(i,m)||_{q(m)} < \infty,
$$
$$
G = \left(\sum_{m=1}^d 1/q(m)  \right)^{-1} \in (0,\infty), \  \ K = \prod_{m=1}^d K(m),
$$
as long as $ \sup_I ||\xi(I)||_G \le C(\vec{q}) \prod_{m=1}^d K(m). $ We receive in this case, 
i.e. if  both  conditions (6.1), (6.2) are satisfyed:
$$
T(Q_d,x) \le C_1 \exp \left( - C_2 (x/K)^G \right),
$$
 and this estimation is not improvable, e.g. for the polynomial $ Q_d = \prod_{m=1}^d \xi(1,m). $
 Thus,  we can assume: 
$$
\sum_{I \in I(d,n)} |b(I)| = \infty,  \ \ \sum_{I \in I(d,n)} b^2(I) < \infty. \eqno(6.3)
$$
(Note that the condition (6.3) is not trivial only if $ n = \infty $).
 We introduce a two measures on the subsets $ I(d,n): $
$$
\mu_1(A) = \sum_{I \in A} |b(I)|, \  \ \mu_2(A) = \sum_{I \in A} b^2(I),
$$
and  introduce for  $ \lambda = const > 0 $ the functions $ \ A(\lambda) = \{I: |b(I)| \le \lambda \}, 
\ B(\lambda) = \{I: |b(I)| > \lambda \},\  a_1(\lambda) = \mu_1(A(\lambda)), \ 
a_2(\lambda) = \sqrt{\mu_2(B(\lambda))}. 
$ \\
{\bf Theorem 13.} {\it For some} $ C_{1,2} = C_{1,2}(b(\cdot)) \in (0, \infty) \ \ \ T(Q_d,x) \le $
$$
 C_1 \inf_{\lambda > 0}  \left( \exp \left(-C_2 \left(x/(a_1(\lambda) K) \right) \right)^G +
\exp  \left(-C_2 \left( x/(a_2(\lambda) K) \right)^{M(d,q)} \right) \right).
$$
{\bf Theorem 14.} {\it If in addition  all the variables} $ \{\xi(i,m) \} $ {\it are independent 
then}   $ T(R_d,x) \le $
$$
 C_1 \inf_{\lambda > 0}  \left( \exp \left(-C_2 \left(x/(a_1(\lambda) K )\right) \right)^G + 
\exp  \left(-C_2 \left( x/(a_2(\lambda) K) \right)^{N_d(q)}  \right) \right). 
$$
{\bf Proof}. We shall assume without loss of generality  $ ||\xi(i,m)||_{q(m)} = 1. $
 Let $ A $ be some subset of $ I(d,n) $ and $ \overline{A} = I(d,n) \setminus A. $  We write: 
$ Q_d = Q_d(1) + Q_d(2), $ where 
$$
Q_d(1) = \sum_{I \in A} b(I) \xi(I), \ \ \ Q_d(2) = \sum_{I \in \overline{A}} b(I)\xi(I).
$$
We get for $ Q_d(1) \ - $ sum from triangule inequality:
$$
||Q_d(1)||_G \le C_3(q) \sum_{I \in A} |b(I)| \ || \xi(I)||_G \le C_4(q) \ \mu_1(A),
$$
hence 
$$
T(Q_d(1),x) \le C_5(q) \exp \left( - C_6 (x/\mu_1(A))^G \right). 
$$
  We obtain for the second sum $ Q_d(2) $ by virtue of theorem 1:
$$
||Q_d(2)||_M \le C_7  \sqrt{\sum_{I \in \overline{A} } b^2(I) } = C_7 
\sqrt{\mu_2( \overline{A})}.
$$
 The proposition of the theorem 13 follows from the elementary inequality 
$$
T(Q_d,x) \le T(Q_d(1),x/2) + T(Q_d(2),x/2)
$$ 
after the minimization over  set $ A; $ it is easily to see that the optimal choosing $ A $ 
has a view $ A = A(\lambda) $ for some $ \lambda \ge 0. $ The proof of  theorem 14 is like to one.\par
{\it Example}. Assume that $  q(1)=q(2)=\ldots = q(d) = q = const \in (0,\infty),$
so that 
$$
T(\xi(i,m),x) \le \exp \left( -x^q \right), \ x \ge 0;
$$
and assume also $ |b(I)| \le C  |I|^{-\alpha}, \ \alpha > d/2; \ $ here
$$
|I| = |(i(1),i(2),\ldots, i(d))| = \sqrt{ \sum_{m=1}^d i^2(m) }.
$$ 
 We deduce from theorem 13 (in the martingale case): if $ \alpha \in (d/2,d), $ then
$$
T(\sum b(I) \xi(I),x) \le C_1(d,q,\alpha) \exp \left( -C_2(d,q,\alpha) x^{ q / (q(d-\alpha) 
+ d) } \right).
$$
 If $ \alpha = d, $ then we conclude
$$
T(Q_d,x) \le C_1 \exp \left(-C_2 (x/\log x)^G \right), \ x \ge 3.
$$
 We conclude in the case  $ \alpha > d $ 
$$
T(Q_d,x) \le C_1 \exp \left(-C_2 x^G \right), \ x \ge 0.
$$
 If $ d = 1 $ and $ q = + \infty, $ i.e. if  $ \exists K \in (0,\infty), \ \forall x > K \ \Rightarrow  
T(\xi(i,m),x) = 0, $  yields a well - known result ( see, e. g. [28], p. 33 - 37.)\par
 Now we shall consider the {\it moment} estimations for $ Q_d $ in our case. Assume that 
for some $ p \ge 2 $
$$
 \sup_i |\xi(i,m)|_p \le 1, 
$$ 
and that  all the variables $ \{\xi(i,m)\} $ are independent. \\
{\bf Theorem 15.} 
$$
|Q_d|_p \le  C(b(\cdot)) \ \inf_{\lambda > 0} \left( a_1(\lambda) + a_2(\lambda) p^d/\log p) \right).
$$
{\bf Proof} is the same as in theorem 13. We obtain  
 for the sum $ Q_d(1) $ using triangular inequality:
$$
|Q_d(1)|_p \le C \sum_{I \in A(\lambda) } |b(I)| \ |\xi(I)|_p \le C \mu_1(A(\lambda)) =
C a_1(\lambda).
$$
 We deduce from  theorem 7:
$$
|Q_d(2)|_p \le C \sqrt{ \sum_{I \in \overline {A(\lambda)} } \ b^2(I) }  \ p^d/\log p
= C a_2(\lambda) \ p^d/\log p.
$$
{\it Example.} Assume again in addition $ |b(I)| \le C |I|^{-\alpha}, $ where 
$ \alpha > d/2. $ Then if $ \alpha \in (d/2,d) \ \Rightarrow $ 
$$
 |Q_d|_p \le C(d,\alpha) p^{2(d-\alpha)} (\log p)^{2(1-\alpha/d)}.
$$
 If $ \alpha = d $ then 

$$
|Q_d|_p \le C(d) \log p; 
$$
 and in  the case $ \alpha > d \ \Rightarrow $
$$
 \sup_{p \ge 2} |Q_d|_p \le C(d,\alpha) < \infty. 
$$
{\bf Theorem 16. } {\it Assume now that all the sequences } $ \{ \xi(i,m), F(i) \}, m=1,2,\ldots,d $ 
{\it are  centered martingale - differences  ("martingale case") such that } 
$$
\max_{m}\sup_{i} |\xi(i,m)|_{d \cdot p} \le 1
$$ 
{\it and} $ b \in B. $ {\it Then }
$$
|Q_d|_p \le C(b(\cdot),d) \inf_{\lambda > 0} \left(a_1(\lambda) + a_2(\lambda) p^d \right).
$$
 
 For example assume again in addition (in the "martingale case") $ b(I) \sim C |I|^{-\alpha},
\exists \ \alpha > d/2. $ It follows from 
theorem 16  (if $ \mu_m(d \cdot p) \le 1, m = 1,2,\ldots,d $ 
for some $ p \ge 2) $ that if $ \alpha \in (d/2,d),$
$$
|Q_d|_p \le C(d,\alpha) \ p^{2(d-\alpha)}.
$$ 
 In  cases $ \alpha = d $ and $ \alpha > d $ we obtain (in our condition $ \mu(d \cdot p) < \infty) $
the same results as in independent case (see Examples in the theorem 15).\par 
  {\it Now we shall investigate the case}   
$$
 \sum_{I \subset  I(d,\infty)} b^2(I) \prod_{m=1}^d \sigma^2(i_m,m) = \infty.
$$
 In particular, the series for $ Q_d $ may divergent. Let us consider now the {\it naturally}
normed multiply  sum 
$$
\theta_n = \sum_{I \subset I(d,n)} \xi(I)/  \left[\sqrt{\sum_{I \subset I(d,n)} \prod_{m=1}^d 
\sigma^2(i_m,m)} \right],
$$
so that $ {\bf E} \theta_n = 0, \ {\bf D} \theta_n = 1. $ \par
 {\bf Lemma 3.} {\it 1. If}
$$
\sup_{i,m} T(\xi(i,m)/\sigma(i,m), \ x) \le \exp \left(-x^q \right), \ q,x \ge 0,
$$ 

{\it then at} $ x > 2 $
$$
\sup_n T(\theta_n,x) \le \exp \left(- C x^{M(d,q)} \right).
$$

{\it 2. If in addition the variables} $ \xi(i,m) $ {\it are independent, then}  
$$
\sup_n T(\theta_n,x) \le \exp \left( - C x^{N_d(q)}  \right).
$$

{\it 3. If }
$$
\sup_{i,m} |\xi(i,m)|_{p \cdot d} /\sigma(i,m) \le 1, 
$$

{\it then}
$$
\sup_n |\theta_n|_p \le C p^d.
$$

{\it 4. If in additional to the 3} $ \{\xi(i,m)\} $ {\it are independent and} 
$$
\sup_{i,m} |\xi(i,m)|_p /\sigma(i,m) \le 1, 
$$
{\it then}

$$
\sup_n |\theta_n|_p \le C p^d/\log p.
$$

{\bf Proof} is very simple. Substituting $ \xi(i,m) = \sigma(i,m) \ \nu(i,m) $ and choosing 
$$
b(I) = \prod_{m=1}^d \sigma(i_m,m)/ \left[\sqrt{ \sum_{I \subset I(d,n)} \prod_{m=1}^d 
\sigma^2(i_m,m)} \right],
$$
we can write 
$$
 \theta_n = \sum_{I \subset I(d,n)}  b(I) \nu(I),\ \ b \in B(d,n).
$$
 We receive  using our estimations (1.5), (1.8), (3.1) and (3.3)  the proposition of lemma 3.\par
 Our result may be considered as some addition to the limit theorem for martingales (see,
for example, [20], p.58.) \par

\vspace{3mm}
{\bf 7. Applications to the weak compactness.} \
 Assume that the multidimensional sequence of coefficients $ \{b(I) \} $ 
dependent on  some parameter $ t; \ t \in V, \ V $ is an arbitrary set:
$ b(I) = b(I,t).$  Suppose  
$$
\sup_{t \in V} \sum_{I \in I(d,\infty)} b^2(I,t) < \infty, \eqno(7.1)
$$
and introduce the following distance between  two arbitrary  points $ t_1, t_2 \in V: $
$$
r_1(t_1,t_2) = \sqrt{\sum_{I \in J(d,\infty)} (b(I,t_1) - b(I,t_2))^2.}
$$
 We shall  consider a random field 
$$
\tau(t) = \sum_{I \in I(d,\infty)} b(I,t) \xi(I),\eqno(7.2)
$$
(series with random coefficients), where $ \{ \xi(i,m) \} $ is again a sequence
of martingale differences so that for some $ q(m) > 0 $

$$
\sup_i \sup_{t \in V} ||\xi(i,m,t)||_{q(m)} < \infty.
$$
{\bf Theorem 17.} {\it Assume that the metric space} $ (V,r_1) $ {\it is full and } 

$$
\int_0^1 H^{1/M(d,q)} (V,r_1,\varepsilon) \ d\varepsilon < \infty,
$$
{\it  where } $ H(V,r_1,\varepsilon) $ {\it is so - called metric entropy of space} $ V $ 
{\it on the distance} $ r_1: $
$$
H(V,r_1,\varepsilon) = \log N(V,r_1,\varepsilon), \ \ N(V,r_1,\varepsilon) =
$$

$$
= \inf_{\{t_i\} }  \{ card \{t_i: \cup_{t_i} \{t: r_1(t,t_i) \le \varepsilon \} = V \}\}.
$$

 {\it Then the series (7.2) convergent uniformly on the } $ t, \ t \in V $ {\it with probability} 
1, $ {\bf P}(\tau(\cdot) \in C(V,r_1) ) = 1, \ ( C(V,r_1) $ {\it denote the space of 
all} $ r_1 \ - $ {\it continuos functions } $ f: V \to R,) $ {\it and for some }
$ C_{10} > 0 $
$$
T(\sup_{t \in V} |\tau(t)|, x) \le \exp \left(-C_{10} x^{M(d,q)} \right), \ x \ge 2. 
$$
{\bf Proof.}  It follows from theorem 1 that:  
$$
\sup_{t \in V} ||\tau(t)||_{M(d,q)} < \infty.
$$
 Further, since  
$$
\tau(t_1) - \tau(t_2) = \sum_{I \in I(d,\infty)} [b(I,t_1) - b(I,t_2)] \xi(I),
$$
we deduce analogously: 
$$
||\tau(t_1) - \tau(t_2)||_{M(d,q)} \le C_{11} \sqrt{ \sum_{I \in I(d,\infty)} 
[b(I,t_1) - b(I,t_2)]^2 } = C_{11} r_1(t_1,t_2).
$$
 Our proposition it follows from the propositions ( [11], p. 195 - 196, 
\ [28], p. 303 - 306.) \par
 Assume now that the coefficients $ b(I) $ are constants, i.e. does not 
dependent on $ \omega \in \Omega, \ t \in V, $ but the martingale - 
differences $ \xi(i,m) $ are separable functions depending still on some  parameter 
 (parameters) $ t, \ t \in V: $
$$
\xi(i,m) = \xi(i,m,t), \ \ t \in V,
$$
i.e. $ \xi(i,m, \cdot) $ are separable random fields. Suppose that for some 
$ q(m) \in (0,\infty] $ 
$$
\sup_{t \in V} \sum_{m = 1}^d \sup_i ||\xi(i,m,t)||_{q(m)} < \infty. 
$$
 Let us introduce the following distance between $ t_1, t_2 \in V: $
$$
r_2(t_1,t_2) = \sum_{m=1}^d \sup_i ||\xi(i,m,t_1) - \xi(i,m,t_2)||_{q(m)} < \infty.
$$
{\bf Theorem 18.} {\it Suppose that the metric space} $ (V,r_2) $ {\it is 
full and that} 
$$
\int_0^1 H^{1/M(d,q)} (V,r_2,\varepsilon) \ d \varepsilon < \infty.
$$
{\it Then the series }
$$
\zeta(t) = \sum_{I \in I(d,\infty) } b(I) \xi(I,t)
$$
{\it convergent uniformly on the parameter} $ t, \ t \in V, $ {\it the random field}
$ \zeta(t) $ {\it belong to the space} $ C(V,r_2) $ {\it with probability 1, and}
$$
T(\sup_{t \in V} |\zeta(t)|,x) \le \exp \left( - C_{12} x^{M(d,q)} \right), 
\ x \ge 1.
$$
{\bf Proof } of  theorem 18 is full analogous to the theorem 17; we 
need only to consider the difference 
$$
\zeta(t_1) - \zeta(t_2) = \sum_{I \in I(d,\infty)} b(I) [\xi(I,t_1) - \xi(I,t_2)].
$$
 and  use  the theorem 1 and identity 
$$
\prod_{m=1}^d y(m) - \prod_{m=1}^d x(m) = \sum_{m=1}^d (y(m) - x(m)) \cdot
\prod_{k \in A(m)} y(k) \cdot \prod_{l \in B(m)} x(l),
$$
where $ A(m), \ B(m) $ are the sets if indexes so that $ A(m) \cap B(m) =
\emptyset, \ card \ A(m) + card \ B(m) = d-1. $ \par
 Assume that the coefficients $ b(I) $ are  functions on  some parameter 
$ \alpha; \ \alpha \in \{\alpha\}, \ b = b_{\alpha}(I), $ and, for instance, 
$$
 \sup_{\alpha} \sum_I b^2_{\alpha}(I) < \infty,
$$

$$
\int_0^1 H^{1/M(d,q)} (V,r_2,\varepsilon) \ d\varepsilon < 
\infty,
$$
then the family of distributions on the space $ \  C(V,r_2) $

$$
\mu_{\alpha} (A) = {\bf P} (\zeta_{\alpha}(\cdot) \in A), 
$$
$ A \ - $ is some Borel subset $ C(V,r_2), $ are weakly compact; here 
$$
\zeta_{\alpha}(t) = \sum_{I \in I(d,\infty)} b_{\alpha}(I) \xi(I,t).
$$

{\bf Remark 5.} Analogous results may be obtained in the terms of Majorizings 
 Measures (see, for example, [28], p. 314 - 318.)\par

{\bf Remark 6.} Probably it is very interesting to  generalize  our results on 
the so - called {\it non - commutative case,} in the spirit of article [40].\par

\vspace{4mm}

{\bf Acknowledgments.} I am wery grateful to prof. M. Lin for useful 
discussions and support of this investigation.\par

\newpage
\begin{center}
  REFERENCES  \\
\end{center}
\vspace{4mm}

{\sc [1] DE LA PE$\hat{N}$A V.H.} (1999)  A general class of exponential 
inequalities for martingales and ratios. Ann. Probab., {\bf 27}, 537 - 554. \\

{\sc [2] KLASS M.J., NOVICKI P.} (1997) Order of magnitude bounds for expectation of $ \Delta_2 $ 
function of non - negative random bilinear forms and $ U \ - $ statistics. Ann. Probab., {\bf 25},
1471 - 1501.  \\

{\sc [3] PESHKIR G., SHIRJAEV A.N.} (1995) The Khintchin inequalities and martingale expansion sphere 
of their action. Russian Math. Surveys, {\bf 50}, 5, 849 - 904.  \\ 

{\sc [4] DZHAPARIDZE H., VAN ZANTEN J.H.} (2001) On Bernstein - type inequalities for 
martingales. Stochastic Process. Appl., {\bf 93},  $N^o$ 1, 109 - 118. \\ 

{\sc [5] HITZENKO P.} (1990) Upper bounds for the $ L_p $ - norms of Martingales.
Probab. Theory Related Fields. {\bf 86}, 225 - 238. \\ 

{\sc [6] PINELIS I.} (1994) Optimum bounds for the distribution of martingales in 
Banach spaces. Ann. Probab., {\bf 22}, 1679 - 1706. \\

{\sc [7] VAN DE GEER S.} (1995)  Exponential inequalities for martingales, with 
applications to maximum likelihood estimation for counting process.
Ann. Statist., {\bf 23}, 1779 - 1801. \\

{\sc [8] AZUMA K.} (1967)  Weighted sums of certain dependent random variables. 
Totohu Math. J., 357 - 367. \\

{\sc [9] LAIB N.} (1999) Exponential - type inequalities for martingale difference 
sequences. Applications to nonparametric regression estimation. Commun.
 Statist. - Theory, Methods, {\bf 28}, 1565 - 1576. \\

{\sc [10] LESIGN  EMM., DALIBOR VOLNY.} (2001) Large deviations for martingales.
Stochastic Process. Appl., {\bf 96}, 143 - 159. \\

{\sc [11] OSTROVSKY E.I.} (1999) Exponential estimations for Random Fields and its 
applications (in Russian). Obninsk, Russia, OINPE. \\

{\sc [12] KOZATCHENKO YU. V., OSTROVSKY E.I.} (1985) The Banach Spaces of 
random Variables of subgaussian Type. Theory Probab. Math. 
Statist., (in Russian). Kiev, KSU, {\bf 32}, 43 - 57. \\

{\sc [13]. IBRAGIMOV R., SHARAKHMETOV SH.} (1998) On an Exact Constant for
the Rosental Inequality. Theory Probab. Appl., v. {\bf 42}, 294 - 302. \\

{\sc [14] GRIGORIEV P.G.} (2001)  Estimates for Norm of Random Polynomials and Their 
Applications. Math. Zametki, v. {\bf 69} $ N^o 6, $ 868 - 872. (2001)\\

{\sc [15] SAULIS L., STATULIAVICIUS W.} (1989) The Limit Theorems on the 
Great Deviations (in Russian). Vilnius, Mokslas. \\

{\sc [16] KOROLYUK V.S., BOROVSKIKH Yu.V.} (1994) Theory of U \ - 
Statistics. Kluwer Verlag, Dodrecht. \\

{\sc [17] RAO M.M., REN Z.M.} (1991) Theory of Orlicz Spaces. Marcel Dekker, 
New York - Basel - Hon Kong. \\

{\sc [18] BOBROV P.B., OSTROVSKY E.I.} (1997) The Adaptive Estimation 
of a Regression, Density and Spectrum. Probability and 
Statistics, Collective Works  of St. Petersburg Branch of 
the Steklov Math. Institute, (in Russian), (2), v. {\bf 244}, 28 - 45. \\ 

{\sc [19] DHARMADHIKARI  S.W., FABIAN V., JOGDEO K.} (1983) Bounds on the Moments of 
Martingales. Ann. Statist., {\bf 11} $ N^o $ 3, 735 - 739. \\

{\sc [20] HALL P., HEYDE C.C.} (1980) Martingale Limit Theory and Applications.
 Academic Press, New York. \\

{\sc [21] GINE E., LATALA R., ZINN J.} (2000) Exponential and Moment Inequalities for 
$ U \ - $ statistics. Hight Dimensional Probability, II - Progress in 
Probability, Birkhauser, 13 - 35. \\

{\sc [22] IBRAGIMOV R., SHARAKHMETOV SH. and CECEN A.} (2001)  Exact Estimation for 
Moments of Random Bilinear Forms. J. Theoret. Probab., 
 {\bf 14}  $ N^o $ 1, 21 - 36. \\

{\sc [23] BULDYGIN V.V., KOZACHENKO YU.V.} (2000)  Metric Characterization of Random Variables 
and Random Processes. AMS, 678, Providence, R.I. \\ 

{\sc [24] RUBIN H., and VITALE R.A.} (1980) Asymptotic distribution of symmetric 
statistics. Ann. Statist., {\bf 8}, 165 - 180. \\

{\sc [25] YU. ZHANG.} (2001) A Martingale Approach in the Studies of Percolation 
Clusters on the $ Z^d \ $ Lattice. J. Theoret. Probab.,
{\bf 14}, $ N^o $ 1, 165 - 187. \\

{\sc [26] CUN - HUI ZHANG.} (2001) Some Moment and Exponential Inequalities for 
$ V \ - $ statistics with bounded Kernels. J. Theoret. Probab., 
 {\bf 14}, $ N^o $ 2, 511 - 525. \\

{\sc [27] ALON N,  SPENCER J.H, ERDOS P.} (1992) The Probabilistic Method. Wiley - 
Interscience Series in Discrete Mathematics and Optimization. New York,
Chichester, Brisbane, Toronto, Singapore. \\

{\sc [28] LEDOUX M., TALAGRAND M.} (1991) Probability in Banach Spaces. Springer 
Verlag, Berlin, Heidelberg, Toronto, Hong - Kong. \\

{\sc [29] GEISS S.} (2001) Contraction Principles for Vector Valued Martingales with 
Respect to Random Variables Having Exponential Tail with Exponent 
$ 2 < \alpha < \infty. $ J. Theoret. Probab., {\bf 14} $ N^o 1, $ 39 - 59.  \\

{\sc [30] ASTASHKIN S.V.} (1999) Multiply Series in Rearrangement Invariant Spaces.
Funct.  Anal. Applic. $ N^o 2, $  {\bf 33,}  141 - 143.  \\

{\sc [31] BURKHOLDER P.} (1998) Sharp Inequalities for Martingales and Stochastic Integrals.
Asterisque, 157/158; 75 - 96.  \\

{\sc [32] DAVIS B.} (1976) On the $L^p$ norms of stochastic Integrals and other Martingales. 
Duke Math. J. {\bf 43,}  697 - 704. \\

{\sc [33] BROWN H., HOBSON D., ROGERS L.G.}  (2001) The Maximum Maximum of a martingale constrained 
 by an intermediate Law. Probab. Theory Related Fields, {\bf 119}, issue 4, 558 - 578. \\  

{\sc [34] CHAO J.- A. and R. - L. LONG.} (1992) Martingale Transforms and Hardy Spaces. Probab. Theory  
Related Fields, {\bf 91},  issue 3, 399 - 404.  \\

{\sc [35] MATSAK I.K., PLICHKO A.N.}  (1988) The Khinchin Inequality for k - multiply Product of 
independent  random  Variables. Mat. Zametki, {\bf 44,} 690 - 694.  \\

{\sc [36] SURGAILIS D.}  (1981) On infinitely divisible self - similar random Fields. Zeitschrift 
Wahrsch.  Theory Verw. Geb., {\bf 58, }  453 - 477.\\

{\sc [37] HOEFDING W.} (1948) A class of Statistics with asymptotically normal Distributions. Ann.
Statist., {\bf 19, } $ N^o 3, $ 293 - 325.  \\

{\sc [38] ENGEL D.} (1982) The multiply stochastic Integral. Memoirs of the American Math. Soc., 
{\bf 38 }  $ N^o 265, $  1 - 82.  \\  

{\sc [39] CARLEN E., KREE P.} (1991) $ L^p $ Estimations on iterated stochastic Integrals. Ann. 
Probab., {\bf 19} $  N^o 1,$  354 - 396.  \\ 

{\sc [40]. MARIUS JUNGE and  QUANHUA XU } (2003). Noncommutative Burkholder/Rosental Inequalities.
Ann. Probab., {\bf 31}, $ N^o 2, $ 948 - 945.\\

\vspace{3mm}

\hspace{40mm} Ben Gurion University, Beer Sheva,  84105. \\

\hspace{40mm} Beer - Sheva,  Ben Gurion street, 2. P.O.Box 61. \\

\hspace{40mm} ISRAEL \\

\hspace{40mm} e - mail: {\bf Galaostr@cs.bgu..ac.il} \\

\hspace{40 mm} Tel. (872) - 08 - 9451613, (872) -08 - 6461606.\\

\hspace{40 mm} Fax: (872) - 08 - 6472873.\\

\end{document}